\documentstyle[amscd,12pt]{amsart}
\numberwithin{equation}{section}
\theoremstyle{plain}
   \newtheorem{thm}{Theorem}[section]
   \newtheorem{lem}[thm]{Lemma}
   \newtheorem{sublem}[thm]{Sublemma}
   \newtheorem{prop}[thm]{Proposition}
   \newtheorem{cor}[thm]{Corollary}
   \newtheorem{conj}[thm]{Conjecture}
   
\theoremstyle{definition}
   \newtheorem{defn}[thm]{Definition}

   \newtheorem{const}[thm]{Construction}
\theoremstyle{remark}
   \newtheorem{rem}[thm]{Remark}
\begin{document}

\author[S. Saito and K. Sato]{Shuji Saito$^{(1)}$ and Kanetomo Sato$^{(2)}$}
\title[Finiteness theorem on zero-cycles]
{A finiteness theorem for zero-cycles over $p$-adic fields}
\date{December 27, 2009}
\thanks{Appendix A was written by Uwe Jannsen.}
\maketitle
\par
\vspace{-5mm}
\begin{center}
{\small
${}^{(1)}$Department of Mathematical Sciences,
University of Tokyo
}
\end{center}
\par
\vspace{-2mm}
\begin{center}
{\scriptsize
8-1 Komaba 3-chome,
Meguro-ku, Tokyo 153-8914, JAPAN
}
\end{center}
\vspace{-2mm}
\begin{center}
{\scriptsize
e-mail address: sshuji@@msb.biglobe.ne.jp
}
\end{center}
\vspace{1mm}
\begin{center}
{\small
${}^{(2)}$Graduate School of Mathematics,
Nagoya University
}
\end{center}
\par
\vspace{-2mm}
\begin{center}
{\scriptsize
Furocho, Chikusa-ku, Nagoya 464-8602, JAPAN
}
\end{center}
\vspace{-2mm}
\begin{center}
{\scriptsize
e-mail address: kanetomo@@math.nagoya-u.ac.jp
}
\end{center}

\tableofcontents
%
%
%
\def\ad{{\mathrm{ad}}}
\def\cd{{\mathrm{cd}}}      
\def\cl{{\mathrm{cl}}}
\def\codim{{\mathrm{codim}}}
\def\et{{\mathrm{\acute{e}t}}} 
\def\ett{{\mathrm{\acute{e}t}}}
\def\ord{{\mathrm{ord}}}      
\def\CH{{\mathrm{C\hspace{-.1pt}H}}}      
\def\Coker{{\mathrm{Coker}}}
\def\Ker{{\mathrm{Ker}}}    
\def\div{{\mathrm{div}}}    
\def\Gal{{\mathrm{Gal}}}    
\def\Br{{\mathrm{Br}}}    
\def\Gys{{\mathrm{Gys}}}
\def\H{{\mathrm{H}}}        
\def\Hom{{\mathrm{Hom}}}    
\def\id{{\mathrm{id}}}     
\def\Image{{\mathrm{Image}}} 
\def\mod{\;{\mathrm{mod}}\;} 
\def\Pic{{\mathrm{Pic}}}    
\def\Proj{{\mathrm{Proj}}}  
\def\Spec{{\mathrm{Spec}}}
\def\trdeg{{\mathrm{trdeg}}}
\def\Tr{{\mathrm{Tr}}}
\def\ch{{\mathrm{ch}}}

%

\def\cF{{\mathcal F}}
\def\cO{{\mathcal O}}
\def\cI{{\mathcal I}}
\def\cC{{\mathcal C}}
\def\cZ{{\mathcal Z}}
\def\cQS{{\mathcal QS}}
\def\cSP{{\mathcal SP}}
\def\cQSP{{\mathcal QSP}}

%
\def\Lam{\varLambda}
\def\lam{\lambda}
\def\vare{\varepsilon}
\def\k{\kappa}
\def\ep{\epsilon}

\def\fm{{\mathfrak{m}}}
\def\fp{{\mathfrak{p}}}
\def\fO{{\mathfrak{O}}}

\def\ra{\rightarrow}
\def\lra{\longrightarrow}
\def\Lra{\Longrightarrow}
\def\la{\leftarrow}
\def\lla{\longleftarrow}
\def\Lla{\Longleftarrow}
\def\da{\downarrow}
\def\hra{\hookrightarrow}
\def\lmt{\longmapsto}
\def\sm{\setminus}
\def\wt#1{\widetilde{#1}}
\def\wh#1{\widehat{#1}}
\def\spt{\sptilde}
\def\ol#1{\overline{#1}}
\def\ul#1{\underline{#1}}
\def\us#1#2{\underset{#1}{#2}}
\def\os#1#2{\overset{#1}{#2}}
\def\lim#1{\us{#1}{\varinjlim}}
\def\sumd#1#2{\underset{x\in {#1}_{#2}}\bigoplus}

\def\bC{{\mathbb C}}
\def\bR{{\mathbb R}}
\def\bZ{{\mathbb Z}}
\def\bQ{{\mathbb Q}}
\def\bG{{\mathbb G}}
\def\bL{{\mathbb L}}
\def\bN{{\mathbb N}}
\def\bP{{\mathbb P}}
\def\bF{{\mathbb F}}
\def\bH{{\mathbb H}}

\def\zl{{\bZ}_\ell}
\def\ql{{\bQ}_\ell}
\def\nz{\bZ/n\bZ}
\def\lnz{\bZ/\ell^n\bZ}
\def\qzl{\ql/\zl}
\def\qlz{\ql/\zl}
\def\Gm{{\Bbb G}_{\hspace{-1pt}\mathrm{m}}}


\def\indlim#1{\underset{#1}{\varinjlim} \ }
\def\projlim#1{\underset{#1}{\varprojlim} \ }


\def\rmapo#1{\overset{#1}{\longrightarrow}}
\def\rmapu#1{\underset{#1}{\longrightarrow}}
\def\lmapu#1{\underset{#1}{\longleftarrow}}
\def\rmapou#1#2{\overset{#1}{\underset{#2}{\longrightarrow}}}
\def\lmapo#1{\overset{#1}{\longleftarrow}}
\def\isom{\overset{\simeq}{\longrightarrow}}

\def\tC{\widetilde{C\,}}
\def\tR{\widetilde{R\,}}
\def\tX{\widetilde{X}}
\def\tS{\widetilde{S}}
\def\tE{\widetilde{E}}
\def\tZ{\widetilde{Z}}

\def\qaq{\quad\text{ and }\quad}
\def\qfor{\quad\text{for }}
\def\qwith{\quad\text{with }}

\def\Ln{\Lam_n}
\def\Linfty{\Lam_\infty}

\def\EX#1#2#3{E^{#3}_{#1,#2}(X)}
\def\EXL#1#2#3{E^{#3}_{#1,#2}(X,\Lam)}
\def\EXLn#1#2#3{E^{#3}_{#1,#2}(X,\Ln)}
\def\EXinfty#1#2#3{E^{#3}_{#1,#2}(X,\Linfty)}


\def\HL#1#2{H_{#1}(#2,\Lam)}
\def\HLn#1#2{H_{#1}(#2,\Ln)}
\def\HLinfty#1#2{H_{#1}(#2,\Linfty)}
\def\tHL#1#2{\widetilde{H}_{#1}(#2,\Lam)}

\def\Het#1{H_{\et}^{#1}}


\def\KCL#1{KC(#1,\Lam)}
\def\KCLn#1{KC(#1,\Lam_n)}
\def\KCLinfty#1{KC(#1,\Lam_\infty)}
\def\KCnz#1{KC(#1,\nz)}
\def\KCqzl#1{KC(#1,\qzl)}


\def\KHL#1#2{KH_{#1}(#2,\Lam)}
\def\KHLn#1#2{KH_{#1}(#2,\Lam_n)}
\def\KHnz#1#2{KH_{#1}(#2,\nz)}
\def\KHlnz#1#2{KH_{#1}(#2,\lnz)}
\def\KHqlz#1#2{KH_{#1}(#2,\qlz)}
\def\KHLinfty#1#2{KH_{#1}(#2,\Lam_\infty)}

\def\edgehom#1{\ep_{#1}}

\def\qss{quasi-semistable}
\def\qsp{$\cQS$-pair}
\def\qspp{$\cQSP$-pair}

\def\IX{{I(X_s)}}

\def\red{{\mathrm{red}}}
\def\Xsred{X_{s,\red}}
\def\tXsred{\tX_{s,\red}}
\def\Ysred{Y_{s,\red}}

\def\PB#1{\bP_B(#1)}
\def\Peta#1{\bP_{\eta}(#1)}
\def\Ps#1{\bP_{s}(#1)}

\def\GB#1{G_B(#1)}
\def\Geta#1{G_{\eta}(#1)}
\def\Gs#1{G_{s}(#1)}
\def\tGB#1{G_B(#1)'}
\def\tGs#1{G_{s}(#1)'}
\def\ttGs#1{G_{s}(#1)''}
\def\ttregGs#1{G_{s}(#1)_{reg}''}
\def\En{E_n}
\def\tEn{\widetilde{E}_n}

\def\hatotzl{\, \widehat{\otimes} \, \zl}

\def\AlbV{Alb_V}
\def\rf{F}
\def\qf{k}
\def\ri{R}
\def\gk{K}

\newpage
\section*{Introduction}

Let $V$ be a smooth projective geometrically integral variety over a field $\qf$. Let $\CH_0(V)$ be the Chow group of zero-cycles on $V$ modulo rational equivalence, and let $A_0(V)\subset \CH_0(V)$ be the subgroup consisting of cycles of degree $0$. Then we have a homomorphism
\begin{equation*}\label{albmap}
\phi_V : A_0(V) \lra \AlbV(\qf)
\end{equation*}
induced by the Albanese map $V \to \AlbV$, where $\AlbV$ denotes the Albanese variety of $V$ and $\AlbV(\qf)$ denotes the group of $\qf$-rational points on $\AlbV$. By Abel's theorem $\phi_V$ is injective if $\dim(V)=1$. It is Mumford \cite{Mu} who first discovered that the situation is rather chaotic in the case of higher dimension.
If $\dim(V)\geq 2$ and $\qf=\bC$, the kernel $\Ker(\phi_V)$ of $\phi_V$ may be too large to be understood by a standard geometric structure.
On the other hand, the situation in the case that $\qf$ is a field of arithmetic nature presents a striking contrast to the case $\qf=\bC$.
If $k$ is finite, $\Ker(\phi_V)$ has been explicitly determined by geometric class field theory \cite{KS} Proposition 9.
One of the most fascinating and challenging conjectures in arithmetic geometry is a conjecture of Bloch and Beilinson that $\Ker(\phi_V)$ is torsion in case $\qf$ is a number field. Combined with the conjecture that $\CH_0(V)$ is finitely generated for $V$ over a number field, it would imply that $\Ker(\phi_V)$ is finite. Very little is known about the conjecture and we have nothing to say about it in this paper. Instead we consider its analogue over $p$-adic fields, more precisely, the following conjecture:
\begin{conj}\label{mainconj}
Let $V$ be a smooth projective geometrically integral variety over a $p$-adic field $\qf$. Then $\Ker(\phi_V)$ is the direct sum of a finite group and a divisible group.
\end{conj}
\noindent
A result of Mattuck \cite{Ma} implies, by aid of \cite{JS1} Lemma 7.8, that $\AlbV(\qf)$ is isomorphic to the direct sum of a finite group and a $p'$-divisible group. Here a $p'$-divisible group means an abelian group which is divisible by any integer prime to $p$. So the above conjecture implies
\begin{conj}\label{mainconjfinite}
Let $V$ be as in Conjecture {\rm\ref{mainconj}}.
Then $A_0(V)$ is the direct sum of a finite group and a $p'$-divisible group.
\end{conj}
\medskip
In essentially equivalent forms, Conjectures \ref{mainconj} and \ref{mainconjfinite} were made by Colliot-Th\'el\`ene already in 1993 (\cite{CT1}) and also in the more recent notes \cite{CT3}, \cite{CT4}. He chose to state most of them as ``questions" or ``guesses". Some finiteness results in this direction have been known for surfaces (\cite{SaSu} Theorem 2.5, see also \cite{CT1} Th\'eor\`eme 2.1) and some special varieties of higher dimension (\cite{RS}, \cite{PS}, \cite{KoSz}, \cite{CT2} Th\'eor\`eme and \S5). What has been known for surfaces is that for any prime $\ell \ne p$ the group $A_0(X)$ is the direct sum of a finite $\ell$-primary group $G_{\ell}$ and an $\ell$-divisible group. But it had not been known in general that $G_{\ell}=0$ for almost all primes $\ell$. That this is indeed the case is now a consequence of Corollary \ref{finitefel} below.
\par
\bigskip
In this paper we will prove the following affirmative result concerning Conjecture \ref{mainconjfinite}, where $\fO_\qf$ denotes the ring of integers in $\qf$:

\begin{thm}[{{\bf Theorem \ref{mainthmfinites9}}}]\label{mainthmfinite}
Let $V$ be a smooth projective geometrically integral variety over a $p$-adic field $\qf$. Let $\rf$ be the residue field of $\fO_{\qf}$. Assume that $V$ has a regular projective flat model $X$ over $\fO_\qf$ on which the reduced subscheme of the divisor $X \otimes_{\fO_\qf}\rf$ has simple normal crossings. Then Conjecture {\rm \ref{mainconjfinite}} holds true.
\end{thm}
\noindent
This result has not been known even in the case that $X \otimes_{\fO_\qf}\rf$ is smooth. See Definition \ref{defqss} below for the definition of divisors with simple normal crossings. Jean-Louis Colliot-Th\'el\`ene notices that a combination of our theorem and one of de Jong's results on alterations (see \cite{dJ} Theorem 8.2) implies the following result, which had been stated as a conjecture in \cite{CT1}.

\begin{cor}\label{finitefel}
Let $V$ be a smooth projective geometrically integral variety over a $p$-adic field $\qf$. Then $A_0(V)$ is $\ell$-divisible for almost all primes $\ell$.
\end{cor}
\noindent
Indeed the cited theorem implies that there exists a finite field extension 
$L$ of $\qf$, a smooth projective variety $W$ over $L$ equipped with a proper generically finite morphism $\pi : W \to V$ such that $W$ has strict semistable reduction over $L$. In particular, it satisfies the assumption of Theorem \ref{mainthmfinite} over $\fO_L$. Hence Corollary \ref{finitefel} follows from Theorem \ref{mainthmfinite} for $W$ and the fact that the composite $A_0(V) \rmapo{\pi^*} A_0(W) \rmapo{\pi_*} A_0(V) $ is the multiplication by a positive integer (cf.\ \cite{Ful} Example 1.7.4).

\begin{cor}\label{finiteRC}
Let the assumptions be as in Theorem {\rm\ref{mainthmfinite}}. Assume further that $V$ is rationally connected in the sense of Koll\'ar, Miyaoka and Mori. Then $A_0(V)$ is the direct sum of a finite group and a $p$-primary torsion group of finite exponent.
\end{cor}
\noindent
It is indeed known that for a smooth projective rationally connected variety $V$ over a field, the group $A_{0}(V)$ is killed by some positive integer (\cite{CT2} Proposition 11). The corollary thus follows from Theorem \ref{mainthmfinite}. In the particular case that $V$ is a smooth compactification of a connected linear algebraic group, Colliot-Th\'el\`ene \cite{CT2} established the conclusion of Corollary \ref{finiteRC} without even assuming the existence of a regular proper model.
\par
\bigskip
For a smooth integral variety $V$ over a field $\qf$, a standard method to study $\CH_0(V)$ is to use the \'etale cycle map (\cite{SGA4.5} Cycle)
$$
\rho_V : \CH_0(V)/n \lra \Het {2d}(V,\nz(d)),
$$
where $d$ denotes $\dim(V)$, $n$ is an integer prime to $\ch(\qf)$ and $\nz(r)=\mu_n^{\otimes r}$ is the Tate twist of the \'etale sheaf $\mu_n$ of $n$-th roots of unity. The \'etale cohomology groups $\Het *(V,\nz(r))$ are finite if $\qf$ is a $p$-adic field. So by the compatibility between $\rho_V$ and $\phi_V$ stated in \cite{R} Appendix, one can easily show that Conjecture \ref{mainconjfinite} follows if $\rho_V$ is injective for any $n \ge 1$ prime to $p$. We do not use this statement in what follows, because the last naive expectation concerning the injectivity of $\rho_V$ was shattered by an example found by Parimala-Suresh \cite{PS}. In this paper we shed new light on the problem by considering the cycle map for schemes over integer rings.
\par
\bigskip
Let $\ri$ be an excellent henselian discrete valuation ring with residue field $\rf$ and fraction field $\qf$. In what follows, we do not assume that $\qf$ has characteristic zero. Put $B:=\Spec(\ri)$ with the closed point $s:=\Spec(\rf)$ and the generic point $\eta:=\Spec(\qf)$. For a scheme $X$ over $B$, put $X_s:=X\times_B s$ and $X_\eta:=X\times_B\eta$, and let $\Xsred$ be the reduced subscheme of $X_s$. Let $\cQSP$ be the category of regular projective flat schemes $X$ over $B$ on which the reduced divisor $\Xsred$ has simple normal crossings. Let $X$ be an object of $\cQSP$, and let $\CH_1(X)$ be the Chow group of one-dimensional cycles on $X$ (cf.\ Definition \ref{def1-1} below).
For an integer $n\geq 1$ prime to $\ch(\rf)$, there is a cycle class map (cf.\ \eqref{purity} and \eqref{cyclemap} below)
\begin{equation*}\label{cyclemapintro}
\rho_X : \CH_1(X)/n \lra \Het {2d} (X,\nz(d))
  \;\; \text{ with } \;\; d:=\dim(X)-1.
\end{equation*}
The main result, Theorem \ref{mainthmfinite} will be deduced from the following result:
\begin{thm}[{{\bf Theorem \ref{injectivity}}}]
\label{injectivityintro}
If $\rf$ is finite or separably closed, then
$\rho_X$ is bijective.
\end{thm}
\medskip
There are three main ingredients in the proof of Theorem \ref{injectivityintro}. The first one is the Bertini theorem over a discrete valuation ring proved in \cite{JS2} Theorem 1. It affirms the existence of a nice hypersurface section $Y \subset X$ over $B$ which belongs to $Ob(\cQSP)$, whose complement $X-Y$ is affine and for which $Y \cup \Xsred$ is a reduced divisor with simple normal crossings on $X$. Such a pair $(X,Y)$ is called an ample \qspp. In this paper we need its refinement, obtained by a joint work of Jannsen and the first author, where it is allowed to impose an additional condition that $Y$ contains a given regular closed subscheme $W \subset X$ with suitable conditions. We include a proof in \S\ref{sect4} below. The second key ingredient is a vanishing theorem of $\Het q(X-Y,\nz(d))$ for an ample \qspp~ $(X,Y)$ (see Theorem \ref{Lcondition} below). A key to its proof is an extension of the base-change theorem of Rapoport and Zink \cite{RZ} Satz 2.19. We remove the assumption on the multiplicity of the components of $X_s$ using the absolute purity theorem of Gabber \cite{Fu}. Recall that we have the affine vanishing theorem of Artin-Gabber (loc.\ cit.\ \S5), which holds for more general affine schemes. Our Theorem \ref{Lcondition} implies that a stronger vanishing holds in our special situation (cf.\ Remark \ref{Lconditionrem}). The last key ingredient of the proof of Theorem \ref{injectivityintro} is an arithmetic invariant associated to a quasi-projective scheme $X$ over $B$. It arises from a homological complex $\KCnz X$ with $(n,\ch(\rf))=1$ of the form
\begin{multline*}\label{KC}
\sumd X 0 H^{0}_{\et}(x,\nz(-1)) \gets \sumd X 1 H^{1}_{\et}(x,\nz)
\gets \cdots \gets \sumd X a H^{a}_{\et}(x,\nz(a-1)) \gets \cdots,
\end{multline*}
where $X_a$ denotes the set of the points on $X$ of dimension $a$, and the most left term is placed in degree $0$. We define the group $\KHnz a X $ as the homology group of $\KCnz X$ in degree $a$ and put
$$
\KHqlz a X:=\indlim n \KHlnz a X\qfor\text{a prime } \ell\not=\ch(\rf).
$$
\par
\medskip
We construct the complex $\KCnz X$ using localization sequences of a certain \'etale homology theory (see \S\S\ref{sect1}--\ref{sect2} below for details). The complex $\KCnz X$ is identified with the case $j=-1$ of the homological complex $C_n^j(X)$ defined by Kato \cite{K}, up to signs of differential operators (cf.\ Definition \ref{defKC} below):
$$
{\small
\begin{CD}
\sumd X 0 H^{j+1}_{\et}(x,\nz(j)) \gets \sumd X 1 H^{j+2}_{\et}(x,\nz(j+1))
\gets \cdots \gets \sumd X a H^{a+j+1}_{\et}(x,\nz(a+j)) \gets \cdots.
\end{CD}
}
$$
In loc.\ cit., he studied the case $j=0$ and conjectures that $C_n^0(X)$ is exact for any $X\in Ob(\cQSP)$ if $\rf$ is finite (loc.\ cit.\ Conjecture 5.1). In this paper we study the case $j=-1$ and propose the following new conjecture.
\begin{conj}[{{\bf Conjecture \ref{KCconj}}}]
\label{KCintro}
Assume that $\rf$ is finite. Then for $X \in Ob(\cQSP)$, there is an isomorphism
$$
\KHqlz a X \simeq
\left.\left\{\gathered
 (\qlz)^{\IX} \\ 
 0 \\
\endgathered\right.\quad
\aligned
&\text{$a=1$}\\
&\text{$a\not=1$},
\endaligned\right.
$$
where $\IX$ denotes the set of irreducible components of $\Xsred$.
\end{conj}
\noindent
The vanishing of $\KHqlz 2 X$ is equivalent to the surjectivity of the cycle map $\rho_X$ (cf.\ Proof of Theorem \ref{artincor}). If $\dim(X)=2$, $\KHqlz 2 X$ is isomorphic to the $\ell$-primary torsion part of the Brauer group of $X$ (cf.\ Lemma \ref{BR}), and its vanishing is a consequence of the base-change theorem for Brauer groups of relative curves (\cite{Gr} III Th\'eor\`eme 3.1, \cite{CTOP} Theorem 1.7\,(c)). See Theorem \ref{artin} below for details. If $\dim(X)=3$, the vanishing of $\KHqlz 3 X$ implies the injectivity of $\rho_X$ (cf.\ Lemmas \ref{longes1} and \ref{MS}). We will prove the following result as a consequence of Theorem \ref{injectivityintro}, where the dimension of $X$ is arbitrary.
\begin{thm}[{{\bf Theorem \ref{mainthmKC}}}]\label{thmKCintro}
Conjecture {\rm\ref{KCintro}} holds in degree $a \le 3$.
\end{thm}
\medskip
We state some other consequences of Theorem \ref{injectivityintro},
  which will be proved in \S\ref{sect9} below.
\begin{cor}[{{\bf Corollary \ref{mainthmcorsepcls9}}}]\label{mainthmcorsepcl}
Let $X$ be an object of $\cQSP$, and let $Y_1,\dots,Y_N$ be the irreducible components of $\Xsred$. Assume that $\rf$ is separably closed. Then for any integer $n>0$ prime to $\ch(\rf)$, we have
$$
 \CH_1(X)/n  \os{\simeq}{\lra} \underset{1\leq i \leq N}{\bigoplus}\nz, \;\;
\alpha \mapsto \langle \alpha,Y_i \rangle_{1\leq i \leq N},
$$
where $\langle \alpha,Y_i \rangle$ denotes the intersection number {\rm(}see {\rm\S\ref{sect9}} below for its definition{\rm)}.
\end{cor}
\noindent
It would be interesting to compare this result with those in \cite{Da} and \cite{Mad}. In both works some nontrivial classes in $A_{0}(V)$ are detected by studying intersection theory with the special fiber of a regular model over the ring of integers. In \cite{Da} the residue field $\rf$ is finite and in \cite{Mad} it is of the shape ${\bC}(\hspace{-1.5pt}(t)\hspace{-1.5pt})$.

\begin{cor}[{{\bf Corollary \ref{goodreductions9}}}]\label{goodreduction}
Assume that $X$ is smooth and projective over $B$ and that $\rf$ is finite or separably closed. Then for any integer $n>0$ prime to $\ch(\rf)$, we have
$$
\begin{CD}
\CH_0(X_s)/n @<{\simeq}<{i^*}< \CH_1(X)/n @>{\simeq}>{j^*}> \CH_0(X_\eta)/n \quad (X_s \, \os{i}{\hra} \, X \, \os{j}{\hookleftarrow} \, X_\eta).
\end{CD}
$$
\end{cor}
\noindent
If $\rf$ is separably closed, the degree map $\CH_{0}(X_{s}) \to \bZ$ induces an isomorphism $\CH_{0}(X_{s})/n$ $\simeq {\bZ}/n$. When $\rf$ is finite, $A_{0}(X_{s})$ is a finite group computed by higher class field theory \cite{KS}. When moreover $X_{s}$ is separably rationally connected, that computation shows $A_{0}(X_{s})=0$. We thus recover the prime-to-$\ch(\rf)$ part of the main result of \cite{KoSz}, which was obtained by deformation theory. One should also compare the above result with the result of \cite{Ko}, where the base is a henselian discrete valuation ring with arbitrary perfect residue field.
\begin{cor}[{{\bf Corollary \ref{affineLs9}}}]\label{affineLintro}
Let $(X,Y)$ be an ample \qspp~ with $\dim(X)=d+1\ge 2$, and assume that $\rf$ is finite {\rm(}resp.\ separably closed{\rm)}. Then for any integer $n>0$ prime to $\ch(\rf)$, the push-forward map
\[ i_* : \CH_1(Y)/n \lra \CH_1(X)/n \quad (i:Y \hookrightarrow X) \]
is bijective for $d \geq 3$ and surjective for $d=2$ {\rm(}resp.\ bijective for $d \ge 2$ and surjective for $d=1${\rm)}.
\end{cor}
\noindent
This result can be viewed as a hard Lefschetz theorem for $\CH_1(X)/n$. Theorem \ref{mainthmfinite} is an easy consequence of it (see \S\ref{sect9} below).
\par
\bigskip
The first author thanks Prof.\ Uwe Jannsen for inspiration on this work obtained from joint works \cite{JS1}, \cite{JS2} and \cite{JS3} on the Kato conjecture and its relation with motivic cohomology of arithmetic schemes. The authors thank him also for his generous permission to include his note on resolution of singularities for embedded curves as the appendix in our paper. They are grateful to Prof.\ Jean-Louis Colliot-Th\'el\`ene for his helpful comments on the first version of this paper and his lucid expository paper \cite{CTb} on the subject of this paper based on his Bourbaki talk. In particular, Corollary \ref{finitefel} was suggested by him. Thanks are also due to Tam\'as Szamuely for his valuable comments.
Finally the authors are grateful to the referee for offering numerous constructive comments to improve the presentation of our paper very much.
\def\dimc{\dim'}
\def\et{}
\newpage
\section{Homology theory and cycle map}\label{sect1}
\stepcounter{thm}
Unless indicated otherwise, all cohomology groups of schemes are taken over the \'etale topology. Throughout this paper, we fix a discrete valuation ring $\ri$ with residue field $\rf$ and quotient field $\qf$. We put $B:=\Spec(\ri)$, $s:=\Spec(\rf)$ and $\eta:=\Spec(\qf)$. Let $Sch_B^{qp}$ be the category of quasi-projective $B$-schemes and $B$-morphisms. For $X \in Ob(Sch_B^{qp})$, we put $X_s:=X\times_B s$ and $X_\eta:=X\times_B\eta$, and write $\Xsred$ for the reduced subscheme of $X_s$. We define $\dimc(X)$ to be the Krull dimension of a compactification of $X$ over $B$, which does not depend on the choice of compactifications. For an integer $q \geq 0$, we define $X_q$ as the set of all points $x$ on $X$ with $\dimc(\overline{\{x\}})=q$, where $\overline{\{x\}}$ is the closure of $x$ in $X$. Note that a point $x$ on $X_{\eta}$ (resp.\ on $X_s$) belongs to $X_q$ if and only if the residue field $\kappa(x)$ satisfies
$$
\trdeg_\qf(\kappa(x)) = q-1 \quad
  \hbox{(resp.}\ \trdeg_\rf(\kappa(x)) = q).
$$
One can easily check that
\begin{equation}\label{dimension}
X_q\cap Y =Y_q\quad\hbox{ for $Y$ locally closed in $X$}.
\end{equation}
Let $\cC$ be the full subcategory of $Sch_B^{qp}$ consisting of objects {\it whose structural morphisms do not factor through $\eta$}. For $X \in Ob(\cC)$, $\dimc(X)$ agrees with the Krull dimension of $X$, and we simply denote $\dimc(X)$ by $\dim(X)$.
\begin{defn}\label{def1-1}
For $X\in Ob(\cC)$ and $q \in \bZ_{\ge 0}$, let $\cZ_q(X)$ be the free abelian group over the set $X_q$. We define the rational equivalence relation on $\cZ_q(X)$ in the same way as in \cite{Ful} \S1.3, and define the group $\CH_q(X)$ as the quotient of $\cZ_q(X)$ by this equivalence relation. We call $\CH_q(X)$, {\it the Chow group of $q$-cycles modulo rational equivalence}. This group is identified with the cokernel of the divisor map defined in loc.\ cit.\ \S1.2
$$
\div : \bigoplus_{x \in X_{q+1}}\; \kappa(x)^{\times} \lra
  \bigoplus_{x \in X_q} \; \bZ.
$$
\end{defn}
\begin{rem}\label{rem1-1}
For a scheme $X \in Ob(\cC)$ whose structural morphism factors through $s \hra B$, the above definition of the Chow groups agrees with that in \cite{Ful} \S1.3. We did not include schemes over $\eta$ in the category $\cC$ to avoid confusions caused by the fact that for $X\in Ob(Sch_B^{qp})$ lying over $\eta$, our notation $\CH_q(X)$ corresponds to the usual Chow group of $(q-1)$-cycles. We will not deal with Chow groups of any schemes which do not belong to $\cC$ until the end of \S\ref{sect8}.
\end{rem}
\par
\medskip
Throughout this paper we fix a prime number $\ell$ different from $\ch(\rf)$ and write
$$\Ln=\lnz \qfor \text{ an integer } n\geq 1 \qaq
\Linfty=\qzl=\indlim m \Ln.$$
Let $\Lam$ denote either $\Ln$ or $\Linfty$. For $X \in Ob(Sch_B^{qp})$ and $q \in \bZ$, we define
\addtocounter{equation}{2}
\begin{equation}\label{homologytheory}
\HL q X:=H^{2-q}(X, Rf^{!}\Lam)
\end{equation}
where $f$ denotes the structural morphism $X \to B$ and $Rf^!$ denotes Deligne's twisted inverse image functor \cite{SGA4} XVIII Th\'eor\`eme 3.1.4. This assignment gives a homology theory on $\cC$ in the sense of \cite{JS1} Definition 2.1. More precisely, we have the following proposition:
\stepcounter{thm}
\begin{prop}\label{lem1-0}
Let $Sch_{B*}^{qp}$ be the category whose objects are the same as $Sch_B^{qp}$ and whose morphisms are proper $B$-morphisms. Then $\HL q {-}$ is a covariant functor from $Sch_{B*}^{qp}$ to the category of abelian groups for each $q \in \bZ$, and satisfies the following {\rm (1)--(3):}
\begin{enumerate}
\item[$(1)$] $\HL q {-}$ is contravariantly functorial with respect to \'etale morphisms in $Sch_B^{qp}$.
\item[$(2)$] If $i : Y \hra X$ is a closed immersion in $Sch_B^{qp}$ with open complement $j : U \hra X$, then there is a long exact sequence
$$
  \dotsb \os{\delta}{\lra} \HL q Y \os{i_*}{\lra} \HL q X \os{j^*}{\lra} \HL q U \os{\delta}{\lra} \HL {q-1} Y \lra \dotsb,
$$
which we call a localization sequence in what follows.
\item[$(3)$] The localization sequence in $(2)$ is functorial with respect to proper morphisms and \'etale morphisms in the following sense. Suppose we are given cartesian squares in $Sch_B^{qp}$
$$
\begin{CD}
   U' @>{j'}>> X' @<{i'}<< Y' \\
  @V{f_U}V{\quad \;\;\;\, \square}V  @V{f}V{\quad \;\;\;\, \square}V  @V{f_Y}VV\\
   U  @>{j}>> X  @<{i}<< \phantom{,}Y,
\end{CD}
$$
where $i$ $($resp.\ $j)$ is a closed $($resp.\ open$)$ immersion and $U=X - Y$. If $f$ is proper, then there is a commutative diagram with exact rows
$$
\begin{CD}
\dotsb @. \; \os{\delta}{\lra} \; @. \HL q {Y'} @. \; \os{i'_*}{\lra} \; @. \HL q {X'} @. \; \os{(j')^*}{\lra} \; @. \HL q {U'} @. \; \os{\delta}{\lra} \; @. \HL {q-1} {Y'} @. \; \lra \; @. \dotsb\phantom{.} \\
  @. @. @V{f_{Y*}}VV @. @V{f_*}VV @. @V{f_{U*}}VV @. @V{f_{Y*}}VV \\
\dotsb @. \; \os{\delta}{\lra} \; @. \HL q Y @. \; \os{i_*}{\lra} \; @. \HL q X @. \; \os{j^*}{\lra} \; @. \HL q U @. \; \os{\delta}{\lra} \; @. \HL {q-1} Y @. \; \lra \; @. \dotsb. 
\end{CD}
$$
If $f$ is \'etale, then there is a similar commutative diagram with vertical arrows $f_Y^*,f^*,f_U^*$ in the opposite direction.
\end{enumerate}
\end{prop}
\begin{pf}
An outline of the proof of these statements is given in \cite{JS1} Example 2.2. The details are straight-forward and left to the reader. We recall here only the definition of the map $g_*: \HL q {Y'} \to \HL q Y$ for a proper morphism $g : Y' \to Y$ in $Sch_B^{qp}$, which will be useful later in Lemma \ref{lem1-1}. Let $\beta$ and $\beta'$ be the structural morphisms $Y \to B$ and $Y' \to B$, respectively. Then $g_*$ is defined as the composite
\stepcounter{thm}
\stepcounter{equation}
\begin{align}
 g :  \HL q {Y'}& =H^{2-q}({Y'},R(\beta')^!\Lam_B)=H^{2-q}(Y,Rg_*R(\beta')^!\Lam_B)
   \notag \\
  &=H^{2-q}(Y,Rg_*Rg^!R\beta^!\Lam_B)
  \os{\ad_g}{\lra} H^{2-q}(Y,R\beta^!\Lam_B) = \HL q Y,
  \label{def:push}
\end{align}
where $\ad_g$ denotes the adjunction map $Rg_*Rg^! \to \id$ obtained by the properness of $g$ (i.e., $Rg_*=Rg_!$) and the adjointness between $Rg_!$ and $Rg^!$.
\end{pf}
\begin{rem}\label{remBO}
For $X \in Ob(\cC)$ whose structural morphism $f$ factors as $X \os{g}{\to} s \os{i}{\hra} B$, we have
$$
\HL q X = H^{2-q}(X, Rf^{!}\Lam_B) = H^{2-q}(X, Rg^!Ri^!\Lam_B) \simeq H^{-q}(X, Rg^!\Lam(-1)_s),
$$
where the last isomorphism follows from the purity for discrete valuation rings \cite{SGA5} I Th\'eor\`eme 5.1. The last group is a special case of the Bloch-Ogus homology \cite{BO} (2.1).
\end{rem}
\par
\medskip
We next provide the following lemma:
\addtocounter{equation}{2}
\begin{lem}\label{lem1-1}
Assume that $X \in Ob(Sch_B^{qp})$ is integral and regular, and put $d:=\dimc(X)-1$. For $m \in \bZ$, let $\Lam(m)$ be the Tate twist. Then for a closed subscheme $Y \subset X$, there is a canonical isomorphism
\begin{equation}\label{purity}
 \tau_{Y,X} : H^{2d-q+2}_{Y}(X,\Lam(d)) \os{\simeq}{\lra} \HL q Y
\end{equation}
which satisfies the following functorial properties$:$
\begin{enumerate}
\item[$(1)$]Let $f:X' \to X$ be a proper morphism $($resp.\ \'etale morphism$)$ in $Sch_B^{qp}$ with $X'$ integral and regular, and let $Y'$ be a closed subscheme of $f^{-1}Y=X' \times_X Y$ $($resp.\ put $Y':=f^{-1}Y)$. Let $g$ be the composite map $Y' \to f^{-1}Y \to Y$. Then there is a commutative diagram
$$
{\small
\begin{CD}
 H^{2d'-q+2}_{Y'}(X',\Lam(d')) @>\tau_{Y',X'}>{\simeq}> \HL q {Y'} \\
  @V{f_*}VV @VV{g_*}V \\
 H^{2d-q+2}_{Y}(X,\Lam(d)) @>\tau_{Y,X}>{\simeq}> \HL q Y
\end{CD}
}
\quad
\left(\text{resp}.\
{\small
\begin{CD}
 H^{2d-q+2}_{Y'}(X',\Lam(d)) @>\tau_{Y',X'}>{\simeq}> \HL q {Y'} \\
  @A{f^*}AA @AA{g^*}A \\
 H^{2d-q+2}_{Y}(X,\Lam(d)) @>\tau_{Y,X}>{\simeq}> \HL q Y
\end{CD}
}
\right),
$$
where $d'$ in the left square denotes $\dimc(X')-1$ and the arrow $f_*$ in the left square is the push-forward map of cohomology groups with supports $($see the proof below$)$.
\item[$(2)$]
Let $i$ be the closed immersion $Y \hra X$ and let $j$ be the open complement $U \hra X$. Then there is a commutative diagram with exact rows
$$
\begin{CD}
\dotsb @. \; \os{\delta}{\lra} \; @. H^{r}_Y(X,\Lam(d)) @. \; \os{i_*}{\lra} \; @. H^{r}(X,\Lam(d)) @. \; \os{j^*}{\lra} \; @. H^{r}(U,\Lam(d)) @. \; \os{\delta}{\lra} \; @. H^{r+1}_Y(X,\Lam(d)) @. \; \lra \; @. \dotsb\phantom{,} \\
  @. @. @V{\tau_{Y,X}}V{\simeq}V @. @V{\tau_{X,X}}V{\simeq}V @. @V{\tau_{U,U}}V{\simeq}V @. @V{\tau_{Y,X}}VV \\
\dotsb @. \; \os{\delta}{\lra} \; @. \HL q Y @. \; \os{i_*}{\lra} \; @. \HL q X @. \; \os{j^*}{\lra} \; @. \HL q U @. \; \os{\delta}{\lra} \; @. \HL {q-1} Y @. \; \lra \; @. \dotsb,
\end{CD}
$$
where we put $r:=2d-q+2$, for simplicity.
\end{enumerate}
\end{lem}
\begin{pf}
Let $D(X_{\ett})$ be the derived category of complexes of \'etale sheaves on $X$, and let $\alpha$ be the structural morphism $X \to B$. We will prove
\begin{quote}
(3) {\it There is a canonical isomorphism $T_{\alpha} : \Lam(d)_X[2d] \os{\simeq}{\lra} R\alpha^!\Lam_B$ in $D(X_{\ett})$.}
\end{quote}
We first finish the proof of the lemma admitting this statement.
We define the isomorphism $\tau_{Y,X}$ in \eqref{purity} as that induced by $T_\alpha$. The right square in (1) is obviously commutative. We obtain the commutative diagram in (2) by applying a standard localization argument to $T_\alpha$ (cf.\ \cite{JSS} (0.4.2)). For a proper morphism $f:X' \to X$ as in (1), we define the morphism
\[ \Tr_f : Rf_*\Lam(d')_X[2d'] \to \Lam(d)_X[2d] \;\; \hbox{ in } \;\; D(X_{\ett}) \] by the commutative square
\[ \begin{CD}
Rf_*\Lam(d')_X[2d'] @>{T_{\alpha'}}>{\simeq}> Rf_*R(\alpha')^!\Lam_B \\
 @V{\Tr_f}VV  @VV{\ad_f}V \\
\Lam(d)_X[2d] @>{T_\alpha}>{\simeq}> R\alpha^!\Lam_B,
\end{CD} \]
where $\alpha'$ denotes the structural morphism $X' \to B$, and $\ad_f$ denotes the adjunction map $Rf_*Rf^! \to \id$ (see $\ad_g$ in \eqref{def:push}). We define the arrow $f_*$ in the left square of (1) as the map induced by $\Tr_f$ and the cobase-change map $H^r_{Y'}(X',-) \to H^r_Y(X,Rf_*(-))$ \cite{SGA4} XVIII (3.1.13.2). In view of the definition of $g_*$ in \eqref{def:push}, $f_*$ makes the left square of (1) commutative. Thus it remains to show (3), which we prove in what follows.
\par
Since $B$ is regular of dimension one, $\alpha$ is flat or otherwise factors as $X \os{\beta}{\to}s \os{i}{\hra} B$. In the flat case, we define $T_\alpha$ as the adjoint of the trace morphism $\Tr_{\alpha} : R\alpha_!\Lam(d)_X[2d] \to \Lam_B$ due to Deligne (\cite{SGA4} XVIII Th\'eor\`eme 2.9). In the latter case, we define $T_\alpha$ as the adjoint of the composite morphism
$$
\begin{CD}
R\alpha_!\Lam(d)_X[2d]=Ri_*R\beta_! \Lam(d)_X[2d] @>{\Tr_{\beta}}>> Ri_* \Lam(-1)_s [-2] @>{\Gys_i}>> \Lam_B,
\end{CD}
$$
where the last arrow is the Gysin map for $i$ (\cite{Fu} \S1.2). We show that $T_\alpha$ is an isomorphism. Fix a locally closed embedding $\gamma : X \hra \bP = \bP^N_B$ into a projective space over $B$. Let $\pi$ be the natural projection $\bP \to B$. Then by \cite{SGA4.5} Cycle, Th\'eor\`eme 2.3.8 (i), $T_{\alpha}$ factors as
$$
\begin{CD}
  \Lam(d)_X[2d] @>{\Gys_{\gamma}}>> R\gamma^!\Lam(N)_{\bP}[2N] @>{T_{\pi}}>> R\gamma^!R\pi^!\Lam_B =R\alpha^!\Lam_B.
\end{CD}
$$
Here the arrow $\Gys_{\gamma}$ denotes the Gysin morphism for $\gamma$ (\cite{Fu} \S1.2), which is an isomorphism by the absolute purity (loc.\ cit.\ Theorem 2.2.1). The second arrow is an isomorphism by \cite{SGA4} XVIII Th\'eor\`eme 3.2.5. Hence $T_\alpha$ is an isomorphism and we obtain the assertion (3). This completes the proof of Lemma \ref{lem1-1}.
\end{pf}
\par
\addtocounter{thm}{1}
\addtocounter{equation}{1}
\begin{cor}\label{cor1-2}
For $X \in Ob(Sch_B^{qp})$, there is a homologically graded spectral sequence \[ E^1_{a,b}(X,\Lam)=\bigoplus_{x\in X_a} H_{\et}^{a-b}(x,\Lam(a-1)) \Longrightarrow \HL {a+b} X,\] which we call the niveau spectral sequence. This spectral sequence is covariant with respect to proper morphisms in $Sch_B^{qp}$ and contravariant with respect to \'etale morphisms in $Sch_B^{qp}$. See Figure {\rm\ref{figure1}} below for a table of $E^1$-terms of this spectral sequence.
\end{cor}
\begin{pf}
By Proposition \ref{lem1-0} and \cite{JS1} Proposition 2.7, there is a spectral sequence \[ E^1_{a,b}(X,\Lam)=\bigoplus_{x\in X_a} \HL {a+b} x \Longrightarrow \HL {a+b} X, \] whose $d^1$-differential maps arise from by definition the connecting homomorphisms in the localization sequences in Proposition \ref{lem1-0}\,(2). For $x \in X$, $\HL q x$ is defined as follows: \begin{equation*} \HL q x := \varinjlim_{U\subset\ol{\{x\}}} \ \HL q U, \end{equation*} where $\ol{\{x\}}$ denotes the closure of $x$ in $X$ and $U$ runs through all non-empty open subsets of $\ol{\{x\}}$. When $x \in X_a$, the limit on the right hand side is isomorphic to the \'etale cohomology group $H_{\et}^{2a-q}(x,\Lam(a-1))$ by $\tau_{U,U}$ in $\eqref{purity}$ for open regular subsets $U \subset \ol{\{x\}}$ and the compatiblity results in Lemma \ref{lem1-1}\,(1). Thus we obtain the desired spectral sequence. The functoriality assertions follow from the corresponding results in \cite{JS1} Proposition 2.7.
\end{pf}
\begin{figure}[htp]
\setlength{\unitlength}{.7mm}
\begin{picture}(100,120)(-30,-35)
      \put(80,80){$E^1_{d+1,d+1}$} 
      \put(40,40){$E^1_{2,2}$} 
      \put(80,40){$E^1_{d+1,2}$} 
      \put(20,20){$E^1_{1,1}$} 
{\scriptsize
      \put(32,25){$d^1_{2,1}$} 
}
      \put(40,20){$E^1_{2,1}$} 
      \put(80,20){$E^1_{d+1,1}$} 
      \put(0,0){$E^1_{0,0}$} 
      \put(20,0){$E^1_{1,0}$} 
      \put(40,0){$E^1_{2,0}$} 
      \put(80,0){$E^1_{d+1,0}$} 
      \put(0,-20){$E^1_{0,-1}$} 
      \put(20,-20){$E^1_{1,-1}$} 
      \put(40,-20){$E^1_{2,-1}$} 
      \put(80,-20){$E^1_{d+1,-1}$} 
      \put(-41,-5){$\HL {-1} X$} %
      \put(-39,15){$\HL 0 X$} %
      \put(-39,35){$\HL 1 X$} %
      \put(-39,55){$\HL 2 X$} %
      \put(-65,75){$E^1_{a,b}=0$ unless $0 \le a \le d+1$ and $b \le a$}
      \multiput(83,57)(0,6){3}{\circle*{1}}
      \multiput(57,57)(6,6){3}{\circle*{1}}
      \multiput(63,42)(3,0){3}{\circle*{1}}
      \multiput(63,22)(3,0){3}{\circle*{1}}
      \multiput(63,2)(3,0){3}{\circle*{1}}
      \multiput(63,-18)(3,0){3}{\circle*{1}}
      \multiput(10,-33)(5,0){15}{\circle*{1}}
\thinlines
      \put(-67,72){\framebox(107,10)}
      \put(59,42){\vector(-1,0){6}}
      \put(79,42){\vector(-1,0){6}}
      \put(39,22){\vector(-1,0){6}}
      \put(59,22){\vector(-1,0){6}}
      \put(79,22){\vector(-1,0){6}}
      \put(19,2){\vector(-1,0){6}}
      \put(39,2){\vector(-1,0){6}}
      \put(59,2){\vector(-1,0){6}}
      \put(79,2){\vector(-1,0){6}}
      \put(19,-18){\vector(-1,0){6}}
      \put(39,-18){\vector(-1,0){6}}
      \put(59,-18){\vector(-1,0){6}}
      \put(79,-18){\vector(-1,0){6}}

      \multiput(-12,-4)(1,-1){12}{\line(0,1){.5}}
      \multiput(-12,16)(1,-1){12}{\line(0,1){.5}}
      \multiput(8,-4)(1,-1){12}{\line(0,1){.5}}
      \multiput(-12,36)(1,-1){32}{\line(0,1){.5}}
      \multiput(28,-4)(1,-1){12}{\line(0,1){.5}}
      \multiput(-12,56)(1,-1){32}{\line(0,1){.5}}
      \multiput(28,16)(1,-1){12}{\line(0,1){.5}}
      \multiput(48,-4)(1,-1){12}{\line(0,1){.5}}
\end{picture}
\caption{Table of $E^1_{a,b}=E^1_{a,b}(X,\Lam)$ ($d:=\dimc(X)-1$)}
\label{figure1}
\end{figure}
\par
The niveau spectral sequence in Corollary \ref{cor1-2} is a fundamental tool in this paper. We often omit the coefficient $\Lam$ of $E^1_{a,b}(X,\Lam)$. We have $\EX ab 1=0$ for any $(a,b)$ with $a < b$ and
\begin{equation}\label{E1chow}
\hspace{-35pt}
\EX 11 2 =\Coker\big(d^1_{2,1} : \bigoplus{}_{x\in X_2}\, \k(x)^\times\otimes\Lam \to \bigoplus{}_{x\in X_1} \, \Lam\big) \simeq \CH_1(X)\otimes \Lam
\hspace{-30pt}
\end{equation}
for $X \in Ob(\cC)$ (see also Remark \ref{rem1-1}), where $\CH_1(X)$ denotes the Chow group of $1$-cycles defined in Definition \ref{def1-1}. The second isomorphism hinges on the fact that the differential map $d^1_{2,1}$ agrees with the negative of the divisor map in Definition \ref{def1-1}, by Lemma \ref{lem1-1}\,(2) and \cite{JSS} Theorem 1.1.1. We define the cycle class mapping
\begin{equation}\label{cyclemap}
\rho_X : \CH_1(X)\otimes\Lam \to \HL 2 X
\end{equation}
to be the edge homomorphism $E^2_{1,1} \to \HL 2 X$ of the niveau spectral sequence. By the functorial properties of the niveau spectral sequences, we obtain the following proposition, which will be frequently used in later sections:
\addtocounter{thm}{2}
\begin{prop}\label{prop1-1}
\begin{enumerate}
\item[$(1)$] The cycle map \eqref{cyclemap} is covariantly functorial with respect to proper morphisms in $\cC$ and contravariantly functorial with respect to \'etale morphisms in $\cC$.
\item[$(2)$] For a closed immersion $i:Y \hookrightarrow X$ in $\cC$ whose open complement $j:U \hookrightarrow X$ belongs to $\cC$, there is a commutative diagram with exact rows
$$
\begin{CD}
\CH_1(Y)\otimes\Lam @>{i_*}>> \CH_1(X)\otimes\Lam @>{j^*}>> \CH_1(U)\otimes\Lam @>>> 0 \\
@VV{\rho_Y}V @VV{\rho_X}V @VV{\rho_U}V \\
\HL 2 Y @>{i_*}>> \HL 2 X @>{j^*}>> \HL 2 U, \\
\end{CD}
$$
where the lower row is the localization exact sequence in Proposition {\rm \ref{lem1-0}\,(2)}.
\end{enumerate}
\end{prop}
\noindent
The exactness of the upper row in the above diagram follows from the same argument as for \cite{Ful} Proposition 1.8. The details are straight-forward and left to the reader.
\begin{rem}\label{rem1-3}
The results in Proposition \ref{lem1-0}\,--\,Corollary \ref{cor1-2} are covariantly functorial in $\Lam$. In particular, for $X \in Ob(\cC)$ and $m,n \ge 1$ there is a commutative diagram with exact rows
$$
\begin{CD}
\CH_1(X)\otimes\Lam_m @>{\times \ell^n}>> \CH_1(X)\otimes\Lam_{m+n} @>{\mod \ell^n}>> \CH_1(X)\otimes\Lam_n @>>> 0 \\ @VV{\rho_X}V @VV{\rho_X}V @VV{\rho_X}V \\
H_2(X,\Lam_m) @>{\times \ell^n}>> H_2(X,\Lam_{m+n}) @>{\mod \ell^n}>> \HLn 2 X.
\end{CD}
$$
where the lower row is a part of the long exact sequence of homology groups associated with the short exact sequence of \'etale sheaves on $B$
\[\begin{CD} 0 @>>> \Lam_m @>{\times \ell^n}>> \Lam_{m+n} @>{\mod \ell^n}>> \Lam_n @>>> 0. \end{CD}\]
\end{rem}
\par
\bigskip
We introduce here the following terminology, in order to state a main result of this paper.
\begin{defn}\label{defqss}
\begin{itemize}
\item[(1)]
For a regular scheme $X$ and a reduced divisor $D$ on $X$, we say that $D$ has {\it simple normal crossings} if it satisfies the following conditions:
\begin{enumerate}
\item[(1a)]
{\it The irreducible components $D_i$ $(i \in I)$ of $D$ are regular.}
\item[(1b)]
{\it For any non-empty subset $J \subset I$, the scheme-theoretical intersection $\bigcap_{i \in J} D_i$ is regular of codimension $\sharp(J)$ in $X$, or otherwise empty.}
\end{enumerate}
\item[(2)]
We call an object $X\in \cC$ {\it \qss}, if it satisfies the following conditions:
\begin{enumerate}
\item[(2a)]
{\it $X$ is pure-dimensional, regular and flat of finite type over $B$.}
\item[(2b)]
{\it The reduced divisor $\Xsred$ on $X$ has simple normal crossings.}
\end{enumerate}
Note that $\Xsred$ is non-empty by the definition of the category $\cC$.
\item[(3)]
We define the category $\cQS$ (resp.\ $\cQSP$) as the full subcategory of $\cC$ consisting of objects which are quasi-semistable (resp.\ quasi-semistable and projective over $B$).
\end{itemize}
\end{defn}
\noindent
The following theorem is one of the main results of this paper:
\begin{thm}[{{\bf Theorem \ref{injectivityintro}}}]
\label{injectivity}
Assume that $R$ is excellent and henselian and that $\rf$ is finite or separably closed. Then for $X \in Ob(\cQSP)$ and $n \ge 1$, the cycle map $\rho_X : \CH_1(X)/\ell^n \to \HLn 2 X$ is bijective.
\end{thm}
\noindent
The proof of Theorem \ref{injectivity} will be completed in \S\ref{sect8}.
We note here the following lemma, where $R$ may not be excellent or henselian and $\rf$ may be arbitrary:
\begin{lem}\label{injlowdim}
If $X\in Ob(\cC)$ has dimension $\leq 2$, then the map $\rho_X : \CH_1(X)/\ell^n \to \HLn 2 X$ is injective for any $n \ge 1$.
\end{lem}
\begin{pf}
This follows from the definition of $\rho_X$ and the fact $\EX ab 1=0$ for $a\geq 3$.
\end{pf}

\newpage

\section{Kato homology}\label{sect2}

The notation remains as in \S\ref{sect1}. In this section, we assume that $R$ is henselian, and introduce an arithmetic invariant associated to each object of $\cC$, which plays a crucial role in the proof of Theorem \ref{injectivity}.

\begin{defn}\label{defKC}
Let $X$ be an object of $\cC$. Let $\Lam$ denote $\Ln$ or $\Linfty$. Let us recall the niveau spectral sequence obtained in Corollary \ref{cor1-2}:
\stepcounter{equation}
\begin{equation}\label{BOss}
E^1_{a,b}(X,\Lam)=\bigoplus_{x\in X_a} H_{\et}^{a-b}(x,\Lam(a-1)) ~~\Longrightarrow ~~
 \HL {a+b} X.
\end{equation}
\begin{itemize}
\item[(1)]
We define the complex $\KCL X$ as the following sequence of $E^1$-terms:
$$
\begin{CD}
  \EXL 0 0 1 @.~\os{d^1_{1,0}}{\lla} ~@. \EXL 1 0  1 @.~\os{d^1_{2,0}}{\lla} ~@. \cdots @.~\os{d^1_{a,0}}{\lla} ~@. \EXL a 0 1 @.~\os{d^1_{a+1,0}}{\lla} ~@. \cdots, \\
  @| @. @| @. @. @. @| \\
\sumd X 0 H^{0}_{\et}(x,\Lam(-1)) @. @. \sumd X 1 H^{1}_{\et}(x,\Lam)
 @. @. @. @. \sumd X a H^{a}_{\et}(x,\Lam(a-1)) @. @.
\end{CD}
$$
where $\EXL a 0 1$ is placed in degree $a$. By Lemma \ref{lem1-1}\,(2) and \cite{JSS} Theorem 1.1.1, the differential maps agree with the negative of the boundary maps of Galois cohomology groups defined in \cite{K} \S1.
In particular, $\KCL X$ is a special case of the Bloch-Ogus-Kato complexes defined in loc.\ cit., up to signs of differential maps.
\item[(2)]
We define the group $\KHL a X$, the Kato homology of $X$, to be $\EXL a 0 2$, which is by definition the homology group of the complex $\KCL X$ in degree $a$.
\end{itemize}
\end{defn}
\par
\medskip
By the functorial properties of the niveau spectral sequences (cf.\ Corollary \ref{cor1-2}), the complex $\KCL X$ is covariantly functorial for proper morphisms in $\cC$ and contravariantly functorial for \'etale morphisms in $\cC$. In particular, for a closed immersion $i:Y \hookrightarrow X$ in $\cC$ whose open complement $j:U \hookrightarrow X$ belongs to $\cC$, there is an exact sequence of complexes (cf.\ \eqref{dimension})
\begin{equation}\label{KCexactsequence}
0 \lra \KCL Y \rmapo{i_*} \KCL X \rmapo{j^*} \KCL U \lra 0.
\end{equation}
By Lemma \ref{E1vanish} below, there is an edge homomorphism
\begin{equation}\label{edgehom}
\edgehom a \;:\; \HLinfty {a} X \lra \KHLinfty a X 
\quad\text{ if $\rf$ is finite},
\end{equation}
\begin{equation}\label{edgehomsepcl}
\edgehom a \;:\; \HLn {a} X \lra \KHLn a X 
\quad\text{ if $\rf$ is separably closed}.
\end{equation}
\addtocounter{thm}{4}
\begin{lem}\label{E1vanish}
Assume that $R$ is henselian. Let $X$ be an object of $\cC$.
\begin{itemize}
\item[(1)]
If $\rf$ is separably closed, then $\EXLn a b 1=0$ for $b\leq -1$.
\item[(2)]
If $\rf$ is finite, then $\EXLn a b 1=0$ for $b\leq -2$ and $\EXinfty a {-1} 1=0$.
\end{itemize}
\end{lem}
\begin{pf}
(1) and the first vanishing of (2) follow from 
the fact $\cd(x)=a+\cd(\rf)$ for $x\in X_a$. The second vanishing of (2) 
follows from the following lemma, which is an analogue of a result of Kahn \cite{Ka} p.\ 61 Proposition 4.
\end{pf}
\begin{lem}\label{E1vanish2}
Assume that $R$ is henselian and that $\rf$ is finite. Let $r\geq 0$ be an integer. Let $L$ be a field which is either finitely generated over $\qf$ with $\trdeg_\qf(L)=r-1$ or finitely generated over $\rf$ with $\trdeg_\rf(L)=r$. Then we have $H^{r+1}(L,\Linfty(n))=0$ for any $n\not=r$. Here for a field $\gk$, $H^*(\gk,-)$ denotes the \'etale cohomology of $\Spec(\gk)$, which agrees with the Galois cohomology of the absolute Galois group $G_\gk$ $($cf.\ \cite{Mi} {\rm III} Example {\rm 1.7)}.
\end{lem}
\begin{pf}
We prove the case that $L$ is over $\qf$. The case that $L$ is over $\rf$ is similar and left to the reader (see also Lemma \ref{claim1} below). If $L$ is a finite extension of $\qf$, then the assertion follows from the Tate duality for Galois cohomology of local fields (cf.\ \cite{Se} II \S5.2 Th\'eor\`eme 2). We proceed the proof by induction on $r=\trdeg_{\qf}(L)+1$. Assume $r \ge 2$. There is a subfield $L_{r-1} \subset L$ which is purely transcendental over $\qf$ and for which $L/L_{r-1}$ is finite. Take a subfield $L_{r-2} \subset L_{r-1}$ which is purely transcendental over $\qf$ of degree $r-2$, and let $\gk$ be the algebraic separable closure of $L_{r-2}$ in $L$. Note that $\trdeg_\gk(L)=1$ and that $\trdeg_{\qf}(\gk)=r-2$. Fix a separable closure $\ol \gk$ of $\gk$. There is a Hochschild-Serre spectral sequence
$$
E_2^{a,b}=H^a(\gk,H^b(L\otimes_\gk \ol \gk,\Linfty(n)))
  \Rightarrow H^{a+b}(L,\Linfty(n)).
$$
Since $\cd(\gk) \le \trdeg_{\qf}(\gk)+2=r$ (loc.\ cit.\ II \S\S4.2--4.3) and $\cd(L\otimes_\gk \ol \gk)=1$ (loc.\ cit.\ II \S\S3.2--3.3), we have $E_2^{a,b}=0$ if $b\geq 2$ or $a\geq r+1$. Hence it suffices to show
$$
E_2^{r,1}=H^r(\gk,H^1(L\otimes_\gk \ol \gk,\Linfty(n)))=0.
$$
By Hilbert's theorem 90 (loc.\ cit.\ II \S1.2 Proposition 1), we have
\begin{align*}
H^r(\gk,H^1(L\otimes_\gk \ol \gk,\Linfty(n)))&=
H^r(\gk,(L\otimes_\gk \ol \gk)^{\times} \otimes \Linfty(n-1))\\
&=\indlim {\gk'\subset\ol \gk}
H^r(\gk,(L\otimes_\gk \gk')^\times \otimes \Linfty(n-1)),
\end{align*}
where $\gk'\subset \ol \gk$ ranges over all finite Galois extensions of $\gk$. See loc.\ cit.\ II \S2.2 Proposition 8 for the last equality. We show that the last group is zero. Fix a finite Galois extension $\gk'$ of $K$. Since $\cd(\gk) \le r$, the corestriction map
$$
H^r(\gk',(L\otimes_\gk \gk')^\times \otimes \Linfty(n-1)) \lra H^r(\gk,(L\otimes_\gk \gk')^\times \otimes \Linfty(n-1))
$$
is surjective. Because $G_{\gk'}$ acts on $(L\otimes_\gk \gk')^\times$ trivially, the natural map induced by cup-product
$$
(L\otimes_\gk \gk')^\times \otimes H^r(\gk', \Linfty(n-1)) \lra H^r(\gk',(L\otimes_\gk \gk')^\times \otimes \Linfty(n-1))
$$
is bijective. Here the left hand side is zero by the induction hypothesis. Thus we obtain Lemmas \ref{E1vanish2} and \ref{E1vanish}.
\end{pf}
\par
\medskip
Assume that $\rf$ is finite and put $\Lam := \Lam_{\infty}$ until we state Conjecture \ref{KCsepcl}. Let $D(\zl)$ be the derived category of complexes of $\zl$-modules. For $X \in Ob(\cQSP)$, we construct a morphism
\addtocounter{equation}{2}
\begin{equation}\label{tracemap}
\KCL X  \lra \Lam^{\IX}[-1]   \;\; \hbox{ in } \;\; D(\zl).
\end{equation}
Here $\IX$ denotes the set of the irreducible components of $\Xsred$ and $\Lam^{\IX}[-1]$ is the homological complex with $\Lam^{\IX}$ placed in degree $1$ and with $0$ in the other degrees. We show that the first boundary map of $\KCL X$ is surjective:
$$
d^1_{1,0} : \underset{y \in X_1}{\bigoplus} H^1(y,\Lam) \to\hspace{-9pt}\to \underset{x\in X_0}{\bigoplus} H^0(x,\Lam(-1)).
$$
Indeed, for a given point $x\in X_0$ there exists a closed point $y$ on $X_{\eta}$ whose closure $\ol {\{ y \}} \subset X$ contains $x$ and is regular at $x$ by Lemma \ref{bertinimoving} below. Because $\ri$ is henselian, $\ol {\{ y \}}$ is the spectrum of a henselian discrete valuation ring. By Hensel's lemma, all roots of unity of order prime to $\ch(\kappa(x))$ in $\kappa(x)$ lift to $\kappa(y)$, and one can easily check the $(y,x)$-factor of $d^1_{1,0}$ is surjective by \cite{JSS} Theorem 1.1.1. Thus $d^1_{1,0}$ is surjective, and we get a natural morphism
\begin{equation}\label{tracemap2}
\KCL X \lra H_1(\KCL X)[-1] \simeq \HL 1 X [-1] \;\; \hbox{ in }\;\; D(\zl),
\end{equation}
where the second isomorphism is given by the map \eqref{edgehom} with $a=1$ (cf.\ Figure 1 after Corollary \ref{cor1-2} and Lemma \ref{E1vanish}). On the other hand, there is a composite map
\begin{equation}\label{H1isom}
\HL 1 X \os{\eqref{purity}}{\simeq} H^{2d+1}_{\et}(X,\Lam(d)) \os{(*)}{\lra} \underset{1\leq j \leq m}{\bigoplus} H^{2d+1}_{\et}(Y_j,\Lam(d)) \isom \Lam^{\IX},
\end{equation}
where $d$ denotes $\dim(X)-1=\dim(X_s)$ and $Y_1,\dotsc, Y_m$ are the irreducible components of $\Xsred$. The arrow $(*)$ is the natural pull-back map, and the last arrow is the trace isomorphisms for projective smooth varieties over $\rf$ (cf.\ \cite{JS1} Lemma 5.3\,(2)). We define the morphism \eqref{tracemap} as the composite of \eqref{tracemap2} and \eqref{H1isom}.
\par
\bigskip
We now propose the following conjectures:
\addtocounter{thm}{3}
\begin{conj}\label{KCsepcl}
If $R$ is excellent and henselian and $\rf$ is separably closed, then $\KCLn X$ is exact for any $X \in Ob(\cQSP)$.
\end{conj}

\begin{conj}[{{\bf Conjecture \ref{KCintro}}}]
\label{KCconj}
If $R$ is excellent and henselian and $\rf$ is finite, then for any $X \in Ob(\cQSP)$, the morphism \eqref{tracemap} is an isomorphism, i.e., we have
\[
\KHLinfty a X \simeq
\left.\left\{\gathered
 \Linfty^{\IX} \\ 
 0 \\
\endgathered\right.\quad
\aligned
&\text{$a=1$}\\
&\text{$a\not=1$.}
\endaligned\right.
\]
\end{conj}
\smallskip
\noindent
Recall that $\KCL X$ is identified with the case $j=-1$ of the Bloch-Ogus-Kato complex
\[ \sumd X 0 H^{j+1}_{\et}(x,\Lam(j)) \gets \sumd X 1 H^{j+2}_{\et}(x,\Lam(j+1)) \gets \cdots \gets \sumd X a H^{a+j+1}_{\et}(x,\Lam(a+j)) \gets \cdots \]
up to signs of differential maps (cf.\ Definition \ref{defKC}). In \cite{K} Conjecture 5.1, he conjectures that this complex with $j=0$ is exact for both $\Lam=\Ln$ and $\Linfty$, if $\rf$ is finite and $X$ is a regular scheme which is proper flat over $B$. It is interesting to ask if Conjecture \ref{KCconj} is valid for $\Ln$. We will prove the following affirmative result on these conjectures using Theorem \ref{injectivity}.

\begin{thm}[{{\bf Theorem \ref{thmKCintro}}}]\label{mainthmKC}
Conjectures {\rm\ref{KCsepcl}} and {\rm\ref{KCconj}} hold in degree $a \le 3$.
\end{thm}
\noindent
The proof of this theorem will be given in \S\ref{sect9} below. In the rest of this section, we prove the following proposition, where $R$ may not be excellent:
\begin{prop}\label{mainthmKClowdim}
Conjectures {\rm\ref{KCsepcl}} and {\rm\ref{KCconj}} hold in degree $a \leq 1$.
\end{prop}
\begin{pf}
Let $\Lam$ denote $\Ln$ if $\rf$ is separably closed and $\Linfty$ if $\rf$ is finite. We have
$$
\KHL a X= \EXL a 0 2 = \EXL a 0 \infty =\HL a X\qfor a\leq 1
$$
by Lemma \ref{E1vanish}. Hence the proposition follows from the following lemma.
\begin{lem}\label{sublemH1isom}
Assume $R$ is henselian, and let $X$ be an object of $\cQSP$.
\begin{enumerate}
\item[(1)]
If $\rf$ is separably closed, then $\HL a X$ is zero for $a \le 1$.
\item[(2)]
If $\rf$ is finite, then the composite map \eqref{H1isom} is bijective and $\HL a X$ is zero for $a \le 0$.
\end{enumerate}
\end{lem}
\medskip
\noindent
{\it Proof of Lemma \ref{sublemH1isom}.}
We prove only the case $a=1$ of (2), i.e., the bijectivity of the arrow $(*)$ in \eqref{H1isom}. The other assertions follow from similar arguments as below and are left to the reader. The map $(*)$ in question factors as follows
$$
H^{2d+1}_{\et}(X,\Lam(d)) \lra H^{2d+1}_{\et}(X_s,\Lam(d)) \lra H^{2d+1}_{\et}(\Xsred,\Lam(d)) \os{(**)}{\lra} \underset{1\leq j \leq m}{\bigoplus} H^{2d+1}_{\et}(Y_j,\Lam(d)),
$$
where all arrows are pull-back maps. Since $R$ is henselian, the first arrow is bijective by the proper base-change theorem \cite{SGA4} XII Corollaire 5.5\,(iii). The second arrow is bijective by the invariance of \'etale topology for the nilpotent closed immersion $\Xsred \hra X_s$ (loc.\ cit.\ VIII Corollaire 1.2). 
We show that the last arrow $(**)$ is bijective. Let $\nu$ be the canonical map $\coprod_{i=1}^m \, Y_j \to \Xsred$, which is finite because each component $Y_j \to \Xsred$ is a closed immersion. Let $\varSigma$ be the singular locus of $\Xsred$. As $\nu$ is an isomorphism over $\Xsred-\varSigma$, the adjunction map $\alpha : \Lam(d) \to \nu_*\nu^*\Lam(d)$ of \'etale sheaves on $\Xsred$ is injective, and the support of $\Coker(\alpha)$ is contained in $\varSigma$. Hence we have a long exact sequence
$$
\begin{CD}
\dotsb @.~ \lra  ~@. H^{2d}_{\et}(\varSigma,\Coker(\alpha)) @.~ \lra  ~@. H^{2d+1}_{\et}(\Xsred,\Lam(d))
  @.~ \os{(**)}{\lra}  ~@. \underset{1\leq j \leq m}{\bigoplus} H^{2d+1}_{\et}(Y_j,\Lam(d)) \\
  @.~ \lra  ~@. H^{2d+1}_{\et}(\varSigma,\Coker(\alpha)) @.~ \lra  ~@. \dotsb, \hspace{60pt} @.
\end{CD}
$$
and the map $(**)$ is bijective by the fact $\cd_{\ell}(\varSigma) \le 2\dim(\varSigma) + \cd_{\ell}(\rf) = 2d-1$ (loc.\ cit.\ X Corollaire 4.3). Thus we obtain the lemma and Proposition \ref{mainthmKClowdim}.
\end{pf}

\newpage

\section{Vanishing theorem}

Let the notation be as in \S\ref{sect1}. For a scheme $X$, $D(X_{\ett})$ denotes the derived category of complexes of \'etale sheaves on $X$.
\begin{defn}\label{deflogpair}
See Definition \ref{defqss}\,(3) for the categories $\cQS$ and $\cQSP$.
\begin{itemize}
\item[(1)]
A \qsp~ is a pair $(X,Y)$ consisting of schemes $X,Y \in Ob(\cQS)$ such that $Y$ is a divisor on $X$ and such that the reduced divisor $\Xsred \cup Y$ on $X$ has simple normal crossings. Note that no component of $Y$ is a component of $\Xsred$, because $Y$ is flat over $B$ by assumption.
\item[(2)]
A \qspp~ is a \qsp~ $(X,Y)$ such that $X,Y\in Ob(\cQSP)$.
For a \qspp~ $(X,Y)$, we call $U:=X-Y$ the complement of $(X,Y)$, and often denote it by $(X,Y;U)$.
\item[(3)]
An {\it ample} \qspp~ is a \qspp~ $(X,Y;U)$ such that $U$ is affine.
\end{itemize}
\end{defn}
\noindent
The aim of this section is to prove the following vanishing result:
\begin{thm}\label{Lcondition}
Let $(X,Y;U)$ be an ample \qspp~ with $\dim(X)=d+1\ge 2$, and assume that $R$ is henselian.
\begin{itemize}
\item[(1)] 
If $\rf$ is separably closed, then $\HLn q U =0$ for $q \leq d+1$.
\item[(2)]
If $\rf$ is finite, then $\HLn q U  =0$ for $q \leq d$.
Furthermore if $d\geq 2$, then $\HLinfty {d+1} U  =0$.
\end{itemize}
\end{thm}

\begin{rem}\label{Lconditionrem}
For a regular scheme $U$ which is flat of relative dimension $d$ over $B$, we have
$$
\HL q U \simeq \Het {2d+2-q}(U,\Lam(d))
$$
by \eqref{purity}. When $U$ is affine and $R$ is henselian, the affine Lefschetz theorems of Artin-Gabber (cf.\ \cite{Fu} \S5) imply that the last group vanishes if $\rf$ is separably closed and $q \leq d$, or if $\rf$ is finite and $q \leq d-1$. Theorem \ref{Lcondition} asserts a stronger vanishing result, assuming that $U$ is the complement of an ample \qspp~ $(X,Y)$. For example, assume $\rf$ is separably closed, and consider the case that $X=\Proj(\ri [T_0,T_1])$ and that $Y \subset X$ is the ample divisor $T_1^m - \pi \cdot T_0^m=0$ ($m \ge 1$), where $\pi$ is a prime element of $\ri$. We compute the natural map $\alpha$ in the exact sequence
$$
 H^2_Y(X,\Ln(1)) \os{\alpha}{\lra} H^2(X,\Ln(1)) \lra H^2(U,\Ln(1))
  \quad (U:=X-Y).
$$
The group $H^2_Y(X,\Ln(1))$ is generated by the divisor class $\cl_X(Y)$ \cite{SGA4.5} Cycle Proposition 2.1.4. Consider the isomorphism $\beta : H^2(X,\Ln(1)) \isom H^2(X_s,\Ln(1)) \isom \Ln$, where the first isomorphism is obtained by the proper base change theorem for \'etale cohomology and the second is the trace map for $X_s \to s$ (cf.\ \eqref{corisom1} below). Then we have $\beta\hspace{1pt}\alpha(\cl_X(Y)) = m \in \Ln$ by Lemma \ref{lemcor1} below. Hence $H^2(U,\Ln(1))$ is non-zero if $\ell$ divides $m$. On the other hand, Theorem \ref{Lcondition} asserts that $H^2(U,\Ln(1))$ is zero if $(X,Y)$ is a \qspp, i.e., $m=1$.
\end{rem}
\par
\medskip
The following lemma is crucial to the proof of Theorem \ref{Lcondition}, which extends the base-change theorem of Rapoport-Zink (\cite{RZ} Satz 2.19) to \qspp s by the absolute purity theorem of Gabber \cite{Fu} Theorem 2.1.1. For a scheme $X \in Ob(\cC)$ and an effective Cartier divisor $Y \subset X$, let $\cl_{X}(Y) \in H^2_Y(X,\Ln(1))$ be the divisor class of $Y$ (\cite{SGA4.5} Cycle \S2.1, \cite{Fu} \S1.1), which induces a morphism $\Ln(-1)_{Y}[-2] \to R\kappa^!\Lam_{n,X}$ in $D(Y_{\ett})$ ($\kappa : Y \hra X$) by the natural identification (cf.\ \cite{SGA4.5} Cat\'egories D\'eriv\'ees II.2.3\,(4))
$$
 H^2_Y(X,\Ln(1)) = \Hom_{D(Y_{\ett})}(\Ln(-1)_Y[-2],Ri^! \Lam_{n,X}).
$$
\begin{lem}\label{RZ}
Let $(X,Y)$ be a \qsp~ over $B$. Put $Z:=\Xsred$ and $W:=\Ysred$. Then$:$
\begin{enumerate}
\item[(1)]
$W$ is an effective Cartier divisor on $Z$ and the divisor class $\cl_Z(W) \in H^2_W(Z,\Ln(1))$ induces an isomorphism $H^{q-2}(W,\Ln(j-1)) \simeq H^q_W(Z,\Ln(j))$ for any $q,j \in \bZ$.
\item[(2)]
Assume further that $R$ is henselian and that $(X,Y;U)$ is a \qspp~ over $B$. Put $V := U_{s,\red}$. Then the pull-back map $H^q(U,\Ln(j)) \to H^q(V,\Ln(j))$ is bijective for any $q,j \in \bZ$.
\end{enumerate}
\end{lem}
\begin{pf}
Put $\Lam:=\Ln$. We first show (1). Since $(X,Y)$ is a \qsp~by assumption, the following square consisting of canonical closed immersions is cartesian:
\addtocounter{equation}{4}
\begin{equation}\label{BCM}
\begin{CD}
W @>{f}>> Z \\
@V{i_Y}V{\quad\;\;\;\;\square}V @VV{i_X}V\\
Y @>{\kappa}>> X. \\
\end{CD}
\end{equation}
In particular, $W$ is an effective Cartier divisor on $Z$ and the divisor class $\cl_Z(W)$ is defined in $H^2_W(Z,\Lam(1))$. Our task is to show that the morphism
$$
   \Gys_f : \Lam(-1)_{W}[-2] \lra Rf^! \Lam_{Z} \;\; \hbox{ in } \;\; D(W_{\ett})
$$
induced by $\cl_Z(W)$ is an isomorphism. For $q \ge 1$, let $Z^{(q)}$ be the disjoint union of $q$-fold intersections of distinct irreducible components of $Z$. We define $W^{(q)}$ ($q \ge 1$) in a similar way. Since $X$ and $Y$ belong to $Ob(\cQS)$, each of $Z^{(q)}$ and $W^{(q)}$ is empty or otherwise regular of pure dimension $\dim(Z)+1-q$ and $\dim(W)+1-q$, respectively. Since $(X,Y)$ is a \qsp, one can easily check
\begin{equation}\label{BCM1}
   W^{(q)} = Z^{(q)} \times_Z W \;\; \hbox{ for any } \, q \ge 1.
\end{equation}
Let $u_q: Z^{(q)} \to Z$ and $v_q : W^{(q)} \to W$ be the canonical finite morphisms. There is an exact sequence of \'etale sheaves on $Z$
\begin{equation}\label{BCM11}
  0 \lra \Lam_Z \os{d^0}{\lra} u_{1*}\Lam_{Z^{(1)}} \os{d^1}{\lra} u_{2*}\Lam_{Z^{(2)}} \os{d^2}{\lra} \dotsb \os{d^{q-1}}{\lra} u_{q*}\Lam_{Z^{(q)}} \os{d^q}{\lra} \dotsb,
\end{equation}
where $d^0$ is the adjunction map and $d^q$ for $q \ge 1$ is an alternating sum of adjunction maps. The signs of these alternating sums are defined as follows. Let $Z_1,\dotsc,Z_m$ be the irreducible components of $Z$, and let $S$ and $T$ be the components of $Z^{(q)}$ and $Z^{(q+1)}$ given by
$$
  S = Z_{i_1} \cap Z_{i_2} \cap \dotsb \cap Z_{i_q} \qaq T = S \cap Z_{i_{q+1}},
$$
where $i_j$'s ($1 \le j \le q+1$) are pairwise distinct and we assumed $i_1 < i_2 < \dotsb < i_q$. Then the sign of the $(S,T)$-factor of $d^q$ is defined as $(-1)^{q-r}$, where $r$ is the integer with $i_r < i_{q+1} < i_{r+1}$. The exact sequence \eqref{BCM11} yields a spectral sequence of sheaves on $W_{\et}$ (cf.\ \cite{Mi} Appendix C\, (g))
$$
  E_1^{a,b}=R^bf^!u_{a*}\Lam_{Z^{(a)}} \Lra R^{a+b-1}f^!\Lam_Z.
$$
We compute the sheaf $E_1^{a,b}$. Since $u_a$ and $v_a$ are finite, we have $R^bf^!u_{a*}=v_{a*}R^b(f^{(a)})^!$ by \eqref{BCM1}, where $f^{(a)}$ denotes the closed immersion $W^{(a)} \hra Z^{(a)}$. Hence we have
$$
 E_1^{a,b} \simeq \begin{cases} v_{a*} \Lam(-1)_{W^{(a)}}  \quad & (\hbox{if }\, b = 2 \, \hbox{ and } \, a \ge 1) \\ 0 & (\hbox{otherwise}) \end{cases}
$$
by the absolute purity \cite{Fu} Theorem 2.1.1, where the isomorphism for $b=2$ arises from the isomorphism $\Lam(-1)_{W^{(a)}}\isom R^2(f^{(a)})^!\Lam_{Z^{(a)}}$ induced by the divisor class $\cl_{Z^{(a)}}(W^{(a)}) \in H^2_{W^{(a)}}(Z^{(a)},\Lam(1))$. In particular, we have $R^qf^!\Lam_Z = E_2^{q-1,2}$ for any $q \in \bZ$. Finally we compute the sheaves $E_2^{*,2}$ by the following diagram of sheaves:
$$
\begin{CD}
0 @.~ \lra ~@. R^2f^!\Lam_Z @.~ \lra ~@. E_1^{1,2} @.~ \os{d_1^{1,2}}{\lra} ~@. E_1^{2,2} @.~ \os{d_1^{2,2}}{\lra} ~@. \dotsb @.~ \os{d_1^{q-1,2}}{\lra} ~@. E_1^{q,2} @.~ \os{d_1^{q,2}}{\lra} ~@. \dotsb \\
  @. @. @A{\Gys_f}AA @. @| @. @| @. @. @. @|\\
0 @.~ \lra ~@. \Lam(-1)_W @.~\os{d^0}{\lra}~@. v_{1*}\Lam(-1)_{W^{(1)}} @.~ \os{d^1}{\lra}~ @. v_{2*}\Lam(-1)_{W^{(2)}} @.~ \os{d^2}{\lra}~ @. \dotsb @.~ \os{d^{q-1}}{\lra}~@. v_{q*}\Lam(-1)_{W^{(q)}}  @. ~ \os{d^q}{\lra} ~@. \dotsb.
\end{CD}
$$
The lower row is an exact sequence analogous to \eqref{BCM11} and the signs of the alternating sum in $d^q$ ($q \ge 1$) are defined by those in \eqref{BCM11} and the isomorphism \eqref{BCM1}. We prove these squares are commutative, which implies Lemma \ref{RZ}\,(1). Indeed, for any $q \ge 1$ the pull-back map $H^2_W(Z,\Lam(1)) \to H^2_{W^{(q)}}(Z^{(q)},\Lam(1))$ sends $\cl_Z(W)$ to $\cl_{Z^{(q)}}(W^{(q)})$ by \eqref{BCM1} and \cite{Fu} Proposition 1.1.3, which shows the commutativity in question. Thus we obtain Lemma \ref{RZ}\,(1).

We prove Lemma \ref{RZ}\,(2). Since $W=Y \times_X Z$ by \eqref{BCM}, the pull-back map $H^2_Y(X,\Lam(1)) \to H^2_W(Z,\Lam(1))$ sends $\cl_X(Y)$ to $\cl_{Z}(W)$ by \cite{Fu} Proposition 1.1.3, which yields a commutative diagram in $D(W_{\ett})$
\begin{equation}\label{BCM2}
\begin{CD}
i_Y^*R\kappa^! \Lam_X @<<< i_Y^*\Lam(-1)_Y[-2]  \\
  @VVV @| \\
 Rf^! \Lam_Z @<{\Gys_f}<< \Lam(-1)_W[-2].
\end{CD}
\end{equation}
Here the left vertical arrow is the base-change morphism associated to the square \eqref{BCM} (cf.\ \cite{SGA4} XVIII 3.1.14.2). The top arrow is the morphism induced by $\cl_X(Y)$, which is an isomorphism by the absolute purity \cite{Fu} Theorem 2.1.1. The bottom arrow is an isomorphism by (1). Hence the left vertical arrow is an isomorphism as well. Since $R$ is henselian, this fact implies Lemma \ref{RZ}\,(2) by similar arguments as for \cite{RZ} Lemma 2.18 and the invariance of \'etale topology for nilpotent closed immersions \cite{SGA4} VIII Corollaire 1.2. The details are straight-forward and left to the reader.
\end{pf}
\par
\medskip
We start the proof of Theorem \ref{Lcondition}. Put $Z:=\Xsred$, $W:=Y_{s,\red}$ and $V:=U_{s,\red}$. By \eqref{purity} and Lemma \ref{RZ}\,(2), we have
\begin{equation}\notag
\HLn q U \simeq H^{2d+2-q}_{\et}(U,\Ln(d))\simeq H^{2d+2-q}_{\et}(V,\Ln(d)).
\end{equation}
By the ampleness assumption, the $F$-variety $V$ is affine of dimension $d$. By the affine Lefschetz theorem (\cite{SGA4} XIV Corollaire 3.2), we have
\begin{equation}\label{homcohU2}
\hspace{-50pt}
H^{2d+2-q}_{\et}(V,\Ln(d)) = 0 \quad
 \begin{cases}
 \hbox{ if $\rf$ is separably closed and $q \le d+1$} \\
 \hbox{ or if $\rf$ is finite and $q \le d$}. \end{cases}
\hspace{-50pt}
\end{equation}
It remains to show $H^{d+1}_{\et}(V,\Linfty(d))$ is zero assuming that $\rf$ is finite and $d\geq 2$. There is a long exact sequence of \'etale cohomology groups
$$
\dotsb \to H^r(V,\Lam_m(d)) \to H^r(V,\Lam_{m+n}(d)) \to H^r(V,\Ln(d)) \to H^{r+1}(V,\Lam_m(d)) \to \dotsb.
$$
Because these groups are finite by \cite{SGA4.5} Th.\ finitude Corollaire 1.10, we obtain the following long exact sequence by standard projective and inductive limit arguments:
$$
\dotsb \to H^r(V,\bZ_{\ell}(d)) \to H^r(V,\bQ_{\ell}(d)) \to H^r(V,\Linfty(d)) \to H^{r+1}(V,\bZ_{\ell}(d)) \to \dotsb.
$$
Since $H^{d+2}_{\et}(V,\bZ_{\ell}(d))=0$ by \eqref{homcohU2}, the map $H^{d+1}_{\et}(V,\bQ_{\ell}(d)) \to H^{d+1}_{\et}(V,\Linfty(d))$ is surjective. Therefore it is enough to show $H^{d+1}_{\et}(V,\bQ_{\ell}(d))$ is zero, which we prove in what follows.
\addtocounter{thm}{5}
\begin{lem}\label{claim1}
For any integer $n \in \bZ$, there is an isomorphism
$$
H^{d+1}_{\et}(V,\ql(n))
 \simeq \Coker\big(1-\varphi_s : H^d(\ol {V},\ql(n)) \to H^d_{\et}(\ol {V},\ql(n))\big).
$$
Here $\varphi_s$ denotes the arithmetic Frobenius operator, i.e., the element of the absolute Galois group $G_\rf$ that corresponds to the Frobenius substitution on a fixed algebraic closure $\ol \rf$ of $F$. For a scheme $X \in Ob(\cC)$ which lies over $s=\Spec(\rf)$, $\ol X$ denotes the scalar extension $X \otimes _\rf \ol \rf$, which we regard as a $G_\rf$-scheme by its natural action on $\ol \rf$.
\end{lem}
\begin{pf*}{\it Proof of Lemma \ref{claim1}}
There is a short exact sequence (cf.\ \cite{CTSS} p.780 (28))
$$
0 \lra H^d(\ol {V},\ql(n))_{G_\rf} \lra H^{d+1}_{\et}(V,\ql(n)) \lra H^{d+1}(\ol {V},\ql(n))^{G_\rf} \lra 0,
$$
where for a $G_\rf$-module $M$, $M^{G_\rf}$ (resp.\ $M_{G_\rf}$) denotes the maximal subgroup (resp.\ maximal quotient group) of $M$ on which $G_\rf$ acts trivially. Since $H^{d+1}(\ol {V},\ql(n))=0$ by \eqref{homcohU2}, we have $H^{d+1}_{\et}(V,\ql(n))\simeq H^d(\ol {V},\ql(n))_{G_\rf}$. The last group is identified with the cokernel of the endomorphism $1-\varphi_s$ on $H^d(\ol {V},\ql(n))$, because $G_\rf$ is topologically generated by $\varphi_s$.
\end{pf*}
By Lemma \ref{claim1}, it is enough to show that the eigenvalues of $\varphi_s$ on $H^d_{\et}(\ol V,\ql(d))$ are different from $1$. By Deligne's theorem on weights (\cite{D} Th\'eor\`eme I) and the assumption that $d\geq 2$, the eigenvalues of $\varphi_s$ on $H^d(\ol Z,\ql(d))$ and $H^{d-1}(\ol W,\ql(d-1))$ are different from $1$. It remains to show that there is a long exact sequence of $G_\rf$-$\ql$-vector spaces
$$
\dotsb \to H^{r-2}(\ol W,\ql(d-1)) \to H^r(\ol Z,\ql(d)) \to H^r(\ol V,\ql(d)) \to H^{r-1}(\ol W,\ql(d-1)) \to \dotsb.
$$
Indeed we have $H^{r-2}(\ol W,\ql(d-1)) \simeq H^r_{\ol W}(\ol Z,\ql(d))$ by Lemma \ref{RZ}\,(1), which is $G_\rf$-equivariant because $\cl_{\ol Z}(\ol W)\in H^2_{\ol W}(\ol Z,\zl(1))$ arises from $\cl_{Z}(W) \in H^2_{W}(Z,\zl(1))$. On the other hand, there is a localization long exact sequence consisting of finite $\Ln$-modules for each $n \ge 1$ \[ \dotsb \to H^{r}_{\ol W}(\ol Z,\Ln(d)) \to H^r(\ol Z,\Ln(d)) \to H^r(\ol V,\Ln(d)) \to H^{r+1}_{\ol W}(\ol Z,\Ln(d)) \to \dotsb. \]Hence we obtain the desired long exact sequence by taking the projective limit of these localization long exact sequences with respect to $n \ge 1$ and then taking $\otimes_{\zl}\ql$. Thus $H^{d+1}(V,\ql(d))$ and $H^{d+1}_{\et}(V,\Linfty(d))$ are zero, and we obtain Theorem \ref{Lcondition}.
\hfill
\qed

\newpage

\section{Bertini theorem over a discrete valuation ring}\label{sect4}

Let the notation be as in \S1. The discrete valuation ring $R$ is not necessarily henselian, and the residue field $\rf$ of $\ri$ is arbitrary in this section.
For a free $\ri$-module $E$ of finite rank, let $$\PB E \lra B=\Spec(\ri)$$
be the associated projective bundle (cf.\ \cite{EGAII} 4.1.1).
Put $\Ps E := \PB E\otimes_\ri \rf$, which is the projective bundle over $\rf$ associated to $E_s:=E \otimes_\ri \rf$. Let $\GB E$ be the set of the invertible (i.e., rank one) $\ri$-submodules $N\subset E$ such that $E/N$ is free. Such $N$ induces a closed immersion
$$
H(N):=\PB {E/N} \hookrightarrow \PB E,
$$
which we call a hyperplane in $\PB E $ corresponding to $N$. There is a specialization map
$$
sp_E  : \GB E \lra \Gs E, \;\; N \mapsto N\otimes_\ri \rf,
$$
where $\Gs E$ is the set of the one-dimensional $\rf$-subspaces of $E\otimes_\ri \rf$. This map is surjective. In terms of hyperplanes in projective bundles, $sp_E$ assigns to a hyperplane $H(N) \subset \PB E$ over $B$, a hyperplane $H_s(N):=H(N) \otimes_\ri \rf \subset \Ps E$ over $s$. We quote here the following fact proved in \cite{JS2} Theorem 1.
\begin{thm}\label{bertini0}
Let $X\in Ob(\cQSP)$ with $\dim(X) \ge 2$.
Then there exist a free $\ri$-module $E$ of finite rank, an embedding $X\hookrightarrow \PB E$ and an invertible $\ri$-module $N\in \GB E$ such that $X\cdot H(N):=X\times_{\PB E} H(N)$ lies in $\cQSP$ and such that $(X,X\cdot H(N))$ is a $\cQSP$-pair.
\end{thm}
\noindent
In this paper, we need the following stronger result, which was obtained also in a joint work of Uwe Jannsen and the first author:
\begin{thm}[{{\bf Jannsen-Saito}}]\label{bertini}
Let $X\in Ob(\cQSP)$ with $\dim(X)=d+1 \ge 2$. Let $Y_1,\dotsc,Y_r$ be the irreducible components of $Y:=\Xsred$, which are smooth of dimension $d$ over $\rf$. For integers $1 \leq a \leq d$ and $1\leq i_i< \cdots< i_a\leq r$, we put $Y_{i_1,\dotsc,i_a} := Y_{i_1}\cap\cdots\cap Y_{i_a}.$ Let $W\subset X$ be a closed subscheme satisfying the following three conditions$:$
\begin{enumerate}
\item[{\rm (i)}]
$W$ is the disjoint union of integral regular schemes $W_1,\dotsc,W_m$.
\item[{\rm (ii)}]
For integers $1 \leq a \leq d$, $1\leq i_i< \cdots< i_a\leq r$ and 
$1\leq \nu\leq m$, 
if $W_\nu \not\subset Y_{i_1,\dotsc,i_a}$, then
$W_\nu \times_X Y_{i_1,\dotsc,i_a}$ is empty or regular of dimension
$< \frac{1}{2} \dim(Y_{i_1,\dotsc,i_a})$.
\item[{\rm (iii)}]
For integers $1 \leq a \leq d$, $1\leq i_i< \cdots< i_a\leq r$ and 
$1\leq \nu\leq m$, 
if $W_\nu \subset Y_{i_1,\dotsc,i_a}$, then
$\dim(W_\nu)< \frac{1}{2} \dim(Y_{i_1,\dotsc,i_a})$.
\end{enumerate}
Then there exist a free $\ri$-module $E$ of finite rank,
an embedding $X\hookrightarrow \PB E$ and
an invertible $\ri$-module $N\in \GB E$ satisfying the conditions:
\begin{itemize}
\item[(1)] 
$X\cdot H(N):=X\times_{\PB E} H(N)$ lies in $\cQSP$ and $(X,X\cdot H(N))$ is a $\cQSP$-pair.
\item[(2)] 
$H(N)$ contains $W$.
\end{itemize}
\end{thm}

\begin{rem}\label{bertinirem}
\begin{itemize}
\item[(1)]
In the original form of Theorem \ref{bertini0}, it is assumed that $\rf$ is infinite. This assumption is removed by virtue of a theorem of Poonen \cite{Po1} Theorem 1.1. By the same token it is another theorem due to him \cite{Po2} Theorem 1.1, that enables us to show Theorem \ref{bertini} in the case that $\rf$ is finite (cf.\ Theorem \ref{AKP} below).
\item[(2)]
Theorem \ref{bertini} will play a crucial role in the proof of Theorem \ref{injectivity} in \S8 below. In fact, in the proof of Theorem \ref{injectivity} for the case that $\rf$ is finite, a standard norm argument enables one to use Theorem \ref{bertini} just in the case that the residue field $\rf$ is infinite.
\item[(3)]
Put $H_s(N):=H(N)\otimes_\ri \rf\subset \Ps E$. Theorem \ref{bertini}\,(i) implies that $H_s(N) \times _{\Ps E} Y_{i_1,\dotsc,i_a}$ is regular for any $1\leq i_i< \cdots< i_a\leq m$ and that it is empty if $\dim(Y_{i_1,\dotsc,i_a})=0$. We also remark that $H_s(N) \times _{\Ps E} Y_{i_1,\dotsc,i_a}$ is not empty if $\dim(Y_{i_1,\dotsc,i_a})\geq 1$, because it is a hyperplane section of $Y_{i_1,\dotsc,i_a}\subset \Ps E$.
\end{itemize}
\end{rem}
\par
\medskip
For the proof of Theorem \ref{bertini}, we need the following fact proved in \cite{JS2} Lemma 1.
\begin{lem}\label{bertinilem}
Let $X\in Ob(\cQSP)$ be endowed with an embedding $X\hookrightarrow \PB E$. Let $H\hookrightarrow \PB E$ be a hyperplane. Assume that $H_s := H \otimes_\ri \rf$ and $Y_{i_1,\dotsc,i_a}$ for any $i_1,\dotsc,i_a$ intersect transversally  on $\Ps E$. Then $X\cdot H \,(=X\times_{\PB E} H)$ lies in $\cQSP$ and $(X,X\cdot H)$ is a $\cQSP$-pair.
\end{lem}
\noindent
Theorem \ref{bertini0} is a rather immediate consequence of this lemma
(cf.\ loc.\ cit.\ \S1).
We deduce Theorem \ref{bertini} from Lemma \ref{bertinilem} by a refined argument.
\par
\medskip
\begin{pf*}{\it Proof of Theorem \ref{bertini}}
Let $W^{f}\subset W$ (resp.\ $W^{nf}$) be the (disjoint) union of those $W_\nu$'s
 which are not contained (resp.\ contained) in $Y$.
By the projectivity assumption on $X$,
 take a free $\ri$-module $E_0$ of finite rank and
  an embedding $i: X\hookrightarrow \PB {E_0}$.
Put $\cO_X(n)=i^*\cO_{\PB {E_0}}(n)$ for an integer $n\geq 1$.
By the Serre vanishing theorem (e.g. \cite{Ha} III Theorem 5.2),
we have
$$
H^1(X,\cO_X(n))=H^1(X,\cO_X(n)\otimes I_X(W^{f}))=0
$$
for a sufficiently large $n>0$,
where $I_X(W^{f})\subset \cO_X$ is the ideal sheaf for $W^{f} \subset X$.
We fix such $n$ in what follows.
The $\ri$-modules
$$
\En := H^0(X,\cO_X(n)) \supset
 \tEn :=H^0(X,\cO_X(n)\otimes_{\cO_X} I_X(W^{f}))
$$
are free of finite rank and we have
$$
\En\otimes_\ri \rf = H^0(X_s,\cO_{X_s}(n)) \qaq
\tEn\otimes_\ri \rf =H^0(X_s,\cO_{X_s}(n)\otimes_{\cO_X} I_{X_s}(W_s^{f})).
$$
Here we have used the flatness of $W^{f}$ over $B$.
The last flatness also implies that
 $I_{X_s}(W_s^{f}) = I_{X}(W^{f}) \otimes_{\cO_X} \cO_{X_s}$,
 so that there is a short exact sequence
$$
\begin{CD}
0\to \cO_X(n)\otimes_{\cO_X} I_X(W^{f})@>{\times \pi}>>
 \cO_X(n)\otimes_{\cO_X} I_X(W^{f})\to 
\cO_{X_s}(n)\otimes_{\cO_{X_s}} I_{X_s}(W_s^{f})\to 0,
\end{CD}
$$
where $\pi$ denotes a prime element of $\ri$.
By this short exact sequence,
 the quotient $\En/\tEn$ is free, and
 moreover, for an invertible $\ri$-submodule
  $N\subset \tEn$ with $\tEn/N$ free,
 the quotient $\En/N$ is again free.
Thus we obtain a commutative diagram
$$
\begin{CD}
\GB \tEn @.~\; \os{\iota}{\hra}~\; @. \GB \En \\
@V{sp_{\tEn}}VV @. @VV{sp_{\En}}V \\
\Gs \tEn @.~\; \os{\iota}{\hra}~\; @. \Gs \En. \\
\end{CD}
$$
The images $\iota(\GB \tEn)$ and $\iota(\Gs \tEn)$
  are identified with the sets
$$
\Phi := \{N \in \GB \En | \; W^{f}\subset H(N)\} \qaq
\Phi_s := \{M \in \Gs \En | \; W^{f}_s \subset H_s(M)\},
$$
respectively.
Here $H(N)$ denotes the hyperplane
  $\PB {E/N} \hookrightarrow \PB E$ and $H_s(M)\subset \Ps {E_n}$
is the hyperplane corresponding to $(E_n\otimes_\ri \rf)/M$.
We have used the embedding $X\hookrightarrow \PB \En$
 associated to $\cO_X(n)$.
Consider the set
\begin{align*}
\Phi_{s,\red}
&:= \{M \in \Gs \En | \; W^{f}_{s,\red}=W^{f}\cap \Xsred \subset H_s(M)\}\\
&\phantom{:}= \{M \in \Gs \En | \; W^{f}\cap Y_{i_1,\dotsc,i_a}\subset H_s(M)
\;\;\text{for any } i_1,\dotsc, i_a \}.
\end{align*}
By definition we have $\Phi_s\subset \Phi_{s,\red}$.
We show the following:
\begin{lem}\label{bertiniclaim}
Put
$$
\cF := \Ker\big(I_{X_s}(\Xsred) \to I_{X_s}(\Xsred)\cdot\cO_{W^{f}_s}\big)
\simeq I_{X_s}(\Xsred)\otimes I_{X_s}(W^{f}_s),
$$
where the tensor products of sheaves are taken over $\cO_{X_s}$.
Replacing $n$ with a sufficiently larger one if necessary,
  assume that
$H^1(X_s,\cO_{X_s}(n)\otimes \cF)=0$.
Then for a given $M\in \Phi_{s,\red}$, there exists an $M'\in \Phi_s$ such that
$$H_s(M)\cap \Xsred = H_s(M')\cap \Xsred \; \subset \Ps \En.$$
\end{lem}
\begin{pf*}{\it Proof of Lemma \ref{bertiniclaim}}
By definition we have
\addtocounter{equation}{5}
\begin{equation}\label{phi}
\Phi_s=\big(\Ker(\tau)-\{0\}\big)/\rf^\times,\quad
\Phi_{s,\red}=\big(\Ker(\sigma \tau)-\{0\}\big)/\rf^\times,
\end{equation}
where $\sigma$ and $\tau$ are as follows:
$$
\begin{CD}
 H^0(X_s,\cO_{X_s}(n)) @>{\tau}>> 
H^0\big(X_s,\cO_{X_s}(n)\otimes\cO_{W^{f}_s}\big)
@>{\sigma}>> H^0\big(X_s,\cO_{X_s}(n)\otimes\cO_{W^{f}_{s,\red}}\big).
\end{CD}
$$
The assumption of the claim implies that the map
$$
H^0\big(X_s,\cO_{X_s}(n)\otimes I_{X_s}(\Xsred)\big)\to 
H^0\big(X_s,\cO_{X_s}(n)\otimes \big(I_{X_s}(\Xsred)\cdot\cO_{W^{f}_s}\big)\big)
=\Ker(\sigma)
$$
is surjective. Hence we have
$$
\Ker(\sigma \tau) = \Ker(\tau) + H^0\big(X_s,\cO_{X_s}(n)\otimes I_{X_s}(\Xsred)\big).
$$
In view of \eqref{phi}, this equality implies the desired assertion.
\end{pf*}
\par
\medskip
We turn to the proof of Theorem \ref{bertini}. Put
\begin{align*}
\Phi_{s,\red} \supset \Psi
&:= \{M \in \Gs \En|\; W^{nf}\cup W^{f}_{s,\red}\subset H_s(M)\}\\
&\phantom{:}= \{M \in \Gs \En|\; W \cap Y_{i_1,\dotsc,i_a}\subset H_s(M)
\qfor \text{ any } i_1,\dotsc, i_a \}
\end{align*}
and let $\Psi_{reg}\subset \Psi$ be the subset of such $M$ that $H_s(M)$ and
$Y_{i_1,\dotsc,i_a}$ intersect transversally on $\Ps {E_n}$ for any $i_1,\dotsc, i_a$.
We use here the following fact, which is due to Altman-Kleiman (\cite{AK} Theorem 7) if $\rf$ is infinite, and due to Poonen (\cite{Po1} Theorem 1.1, \cite{Po2} Theorem 1.1) if $\rf$ is finite.
\addtocounter{thm}{1}
\begin{thm}[{{\bf Altman-Kleiman/Poonen}}]\label{AKP}
Let $Y$ be a smooth projective scheme over $\rf$
 and let $Z\subset Y$ be a smooth closed subscheme.
Assume $\dim(Z)<\frac{1}{2} \dim(Y)$. Then for any integer $n \ge 0$ large enough, there exists a smooth hypersurface section $H\subset X$ of degree $n$ containing $Z$.
\end{thm}
\noindent
By this theorem and the assumptions of Theorem \ref{bertini},
  we may suppose that $\Psi_{reg}$ is non-empty,
    replacing $n$ again with a sufficiently larger one if necessarily.
By the surjectivity of $sp_{\tEn}$ and by Lemma \ref{bertiniclaim},
 there exists an $N\in \Phi \subset \GB\En$ satisfying
  $sp_{\En}(N)\in \Psi_{reg}$.
Now Theorem \ref{bertini} follows from Lemma \ref{bertinilem} 
  for the embedding $X\hookrightarrow \GB \En$.
\end{pf*}
\bigskip

\newpage

\section{Surjectivity of cycle map}

Let the notation be as in \S\ref{sect1}. In this section, we prove the following:
\begin{thm}\label{artin}
Assume that $R$ is henselian and that $\rf$ is finite or separably closed. Let $X$ be a scheme over $\Spec(\ri)$ which belongs to $\cQSP$. Then the cycle map defined in \eqref{cyclemap}
$$
\rho_X : \CH_1(X)\otimes\Lam \to \HL 2 X
$$
is surjective for both $\Lam=\Ln$ and $\Linfty$.
\end{thm}
\begin{pf}
It suffices to prove the case that $\Lam=\Ln$. Without loss of generality, we may assume that $X$ is integral. The case $\dim(X)=1$ follows from the definition of $\rho_X$. We prove the case $\dim(X)=2$. Since $\dim(X)=2$, the niveau spectral sequence \eqref{BOss} gives rise to an exact sequence
$$
0 \lra \CH_1(X)/\ell^n \os{\rho_X}{\lra} \HLn 2 X \lra \KHLn 2 X \lra 0.
$$
We will prove the following:
\begin{lem}\label{BR}
If $\dim(X)=2$, then $\KHLn 2 X$ is canonically isomorphic to $\Br(X)[\ell^n]$, where for a scheme $Z$, $\Br(Z)$ denotes the Grothendieck-Brauer group $\Het 2(Z,\Gm)$.
\end{lem}
\noindent
Since $R$ is henselian and $F$ is finite or separably closed, we have $\Br(X)\{ \ell \}=0$ by \cite{CTOP} Corollaries 1.10\,(b) and 1.11\,(b), which is a consequence the base-change theorem for Brauer groups of relative curves \cite{Gr} III Th\'eor\`eme 3.1, \cite{CTOP} Theorem 1.7\,(c). Hence Theorem \ref{artin} for the case $\dim(X)=2$ follows from this lemma.
\begin{pf*}{Proof of Lemma \ref{BR}}
Let $y$ be the generic point of $X$. Since $\dim(X)=2$, we have
\begin{align*}
\KHLn 2 X & = \Ker\big(d^1_{2,0} : H^2(y,\Ln(1)) \to {\bigoplus}_{x\in X_1}\, H^1(x,\Ln)\big)\\
 & = {\bigcap}_{x \in X_1}\,\Ker\big(d_x : H^2(y,\Ln(1)) \to H^1(x,\Ln)\big)\\
 & = {\bigcap}_{x \in X_1}\,\Br(\cO_{X,x})[\ell^n] = \Br(X)[\ell^n],
\end{align*}
where $d^1_{2,0}$ denotes the differential map of the niveau spectral sequence \eqref{BOss}, and $d_x$ ($x \in X_1$) denotes the $x$-factor of $d^1_{2,0}$. The last equality follows from the assumption that $X$ is regular of dimension $2$ and theorems of Auslander-Goldman and Grothendieck \cite{Gr} II Proposition 2.3, Corollary 2.2. We show the third equality. Fix a point $x \in X^1$, and put $Z:=\Spec(\cO_{X,x})$. By Lemma \ref{lem1-1}\,(2), $d_x$ is identified with the arrow $\delta^2$ in the localization exact sequence of \'etale cohomology groups
$$
 H^1(y,\Ln(1)) \os{\delta^1}{\to} H^2_x(Z,\Ln(1)) \to H^2(Z,\Ln(1)) \to H^2(y,\Ln(1)) \os{\delta^2}{\to} H^3_x(Z,\Ln(1)),
$$
where we have $H^1(y,\Ln(1)) \simeq \kappa(y)^{\times}/\ell^n$ by Hilbert's theorem 90, and $H^2_x(Z,\Ln(1)) \simeq \Ln$ by \cite{SGA4.5} Cycle Proposition 2.1.4. By \cite{JSS} Theorem 1.1.1, $\delta^1$ agrees with the negative of the divisor map $\kappa(y)^{\times}/\ell^n \to \Ln$. Because $\cO_{X,x}$ is a principal ideal domain, $\delta^1$ is surjective and we have $H^2(Z,\Ln(1))=\Ker(\delta^2)$. On the other hand, the short exact sequence $0 \to \Ln(1) \to \Gm \os{\times \ell^n}{\to} \Gm \to 0$ on $Z_{\ett}$ gives rise to a short exact sequence
$$
  0 \lra \Pic(\cO_{X,x})/\ell^n \lra H^2(Z,\Ln(1)) \lra \Br(\cO_{X,x})[\ell^n] \lra 0,
$$
and $\Pic(\cO_{X,x})$ is zero because $\cO_{X,x}$ is local. Thus we have $\Ker(d_x)=\Ker(\delta^2)=H^2(Z,\Ln(1))= \Br(\cO_{X,x})[\ell^n]$, which completes the proof of the lemma.
\end{pf*}

Finally we prove the case that $\dim(X) \geq 3$ by induction on $\dim(X)$. By Theorem \ref{bertini} we take an ample divisor $Y\subset X$ such that $(X,Y)$ is an ample \qspp, and put $U:=X-Y$. By Proposition \ref{prop1-1}\,(2), there is a commutative diagram whose lower row is exact
$$
\begin{CD}
\CH_1(Y)/\ell^n @>>> \CH_1(X)/\ell^n  \\
@V{\rho_Y}VV @V{\rho_X}VV \\
\HLn 2 Y @>>> \HLn 2 X @>>> \HLn 2 U.
\end{CD}
$$
By the induction hypothesis, $\rho_Y$ is surjective. Since $R$ is henselian, we have $\HLn 2 U=0$ by Theorem \ref{Lcondition}, which implies the surjectivity of $\rho_X$. Thus we obtain Theorem \ref{artin}.
\end{pf}
\begin{thm}\label{artincor}
Let the assumptions be as in Theorem {\rm\ref{artin}}.
Let $\Lam$ denote $\Ln$ if $\rf$ is separably closed and $\Linfty$ if
$\rf$ is finite. Then $\KHL 2 X=0$.
\end{thm}
\begin{pf}
The niveau spectral sequence \eqref{BOss} and Lemma \ref{E1vanish} yield isomorphisms
$$\KHL 2 X= \EX 20 2 = \EX 20 \infty= \Coker\big(\rho_X : \CH_1(X)\otimes\Lam \to \HL 2 X \big).
$$
Hence Theorem \ref{artincor} follows from Theorem \ref{artin}.
\end{pf}

\newpage

\section{Blow-up formula}\label{sect6}
Let the notation be as in \S\ref{sect1}. We prove here some lemmas which explain how homology groups behave under blow-ups at closed points. The goal of this section is to prove Lemma \ref{blow-up3} below, which will be used in the proof of Theorem \ref{injectivity}. We do not assume $R$ is henselian in this section. The residue field $\rf$ of $\ri$ is arbitrary in the following lemma.
\begin{lem}\label{blow-up1}
Let $X\in Ob(\cC)$ be regular of dimension $d+1$. Fix a closed point $x$ on $X_s$ and let $\pi:\tX \to X$ be the blow-up of $X$
at $x$ and let $\iota: E \hookrightarrow \tX$ be the exceptional divisor.
For any integer $q$ there is a split short exact sequence
$$
0 \lra \tHL q E \os{\iota_*}{\lra} \HL q \tX \os{\pi_*}{\lra} \HL q X \lra 0.
$$
Here $\tHL q E$ is defined as
$$
\tHL q E :=\Ker(\pi_* : \HL q E \lra \HL q x),
$$
which agrees with $\HL q E$ for $q \ge 1$ because $\HL q x$ is zero for $q \ge 1$ {\rm(}see Remark {\rm\ref{remBO}} or \eqref{purity}{\rm)}.
\end{lem}
\begin{pf}
By Proposition \ref{lem1-0}\,(3), there is a commutative diagram with exact rows
$$
\begin{CD}
\dotsb @.~ \lra ~ @. \HL q E @.~ \lra ~ @. \HL q \tX @.~ \lra ~ @. \HL q {\tX-E} @.~ \lra ~ @. \HL {q-1} E @.~ \lra ~ @. \dotsb \\
 @. @. @VVV @. @V{\pi_*}VV @. @| @. @VVV \\
\dotsb @.~ \lra ~ @. \HL q x @.~ \lra ~ @. \HL q X @.~ \lra ~ @. \HL q {X-x} @.~ \lra ~ @. \HL {q-1} x @.~ \lra ~ @. \dotsb.
\end{CD}
$$
It is enough to show that $\pi_*: \HLn q \tX \to \HLn q X$ is split surjective. We define the map $\pi^*: \HLn q X \to \HLn q \tX$ as follows:
$$
\pi^*: \HLn q X \os{\eqref{purity}}{\simeq} H^{2d+2-q}(X,\Ln(d)) \os{\pi^*}{\lra} H^{2d+2-q}(\tX,\Ln(d)) \os{\eqref{purity}}{\simeq} \HLn q \tX.
$$
We show $\pi_*\pi^*=1$ as endomorphisms on $\HL q X$. By Lemma \ref{lem1-1}\,(1), $\pi_*$ is identified with the push-forward map of \'etale cohomology groups. By the projection formula, i.e., the commutative diagram of cup-product pairings (cf.\ Proof of \cite{Mi} VI Corollary 6.4\,(a))
\begin{equation}\notag
\begin{CD}
H^0(\tX,\Ln) @.~ \times ~@. H^{2d+2-q}(\tX,\Ln(d)) @>{\cup}>> H^{2d+2-q}(\tX,\Ln(d)) \\
  @V{\pi_*}VV @. @A{\pi^*}AA @V{\pi_*}VV \\
H^0(X,\Ln) @. \times @. H^{2d+2-q}(X,\Ln(d)) @>{\cup}>> H^{2d+2-q}(X,\Ln(d)),
\end{CD}
\end{equation}
we have only to check the left vertical $\pi_*$ maps $1$ to $1$, which follows from the fact that $\pi$ is an isomorphism outside of $E$. Thus we obtain the lemma.
\end{pf}
\par
Let $\pi: \tX \to X$ and $\iota : E \hra \tX$ be as in Lemma \ref{blow-up1}. Let $W\subset \tX$ be an integral closed subscheme of dimension one. For $\Lam=\Ln$, let $$\Lam \langle W \rangle \subset \HL 2 \tX$$ be the subgroup generated by $\rho_X(W) \in \HL 2 \tX$. For $\Lam=\Linfty$, we put
$$
\Lam \langle W \rangle:=\indlim {n \in {\mathbb N}} \Ln\langle W \rangle \subset \HLinfty 2 \tX.
$$
See Remark \ref{rem1-3} for the compatibility of cycle class maps with respect to the multiplication by $\ell$ from $\Lam_n$ to $\Lam_{n+1}$.
Recall that $E\simeq \bP^d_{x}$. For an integral closed subscheme $W \subset E$ of dimension one, we define $\deg_{E/x}(W)\in \bZ$ to be the degree over $x$ (cf.\ \cite{Ful} Definition 1.4) of the intersection of $W$ and a hyperplane $H\subset E$ with $W\not\subset H$. Under these settings, we prove:

\begin{cor}\label{blow-up2}
Assume $\cd_{\ell}(\rf) \le 1$\,{\rm(}e.g., $F$ is finite or separably closed{\rm)}. Let $W\subset E$ be a closed integral subscheme of dimension one, and assume either that $\Lam =\Linfty$ or that $\deg_{E/x}(W)$ is prime to $\ell$. 
Then we have
\[ \Ker\big(\pi_* : \HL 2 \tX \to \HL 2 X\big) =  \Lam \langle W \rangle. \]
\end{cor}
\begin{pf}
Lemma \ref{blow-up1} yields an exact sequence
\[ 0 \lra \HL 2 E \os{\iota_*}{\lra} \HL 2 \tX \os{\pi_*}{\lra} \HL 2 X \lra 0. \]
Since $E\simeq \bP^d_{x}$, $\CH_1(E)$ is generated by a line $L \subset E$ and we have
\[ \deg_{E/x} : \CH_1(E) \os{\simeq}{\lra} \bZ. \]
We show that the map $\rho_E : \CH_1(E)\otimes \Lam \to \HL 2 E$ is bijective. Indeed,
we have \[ \HLn 2 E \us{\eqref{purity}}\simeq H^{2d-2}(E,\Ln(d-1)) \us{\cd_{\ell}(\rf) \le 1}\simeq \Ln \] and $\HLn 2 E$ is generated by $\rho_E(L)$, which implies the bijectivity of $\rho_E$ for both $\Lam=\Ln$ and $\Linfty$. Corollary \ref{blow-up2} follows easily from these facts and Proposition \ref{prop1-1}\,(1) for $E \hra \tX$.
\end{pf}
\bigskip
\noindent
{\bf Setting 6.3.}
We consider the following situation:
\addtocounter{equation}{2}
\begin{equation*}\label{sucbu}
\begin{CD}
X=X^{(0)} @<{\pi_1}<< X^{(1)}  @<{\pi_2}<< X^{(2)} @<{\pi_3}<< \dotsb @<{\pi_N}<< X^{(N)} =\tX.
\end{CD}
\end{equation*}
Here $X\in Ob(\cC)$ is regular of dimension $d+1$, and $\pi_i:X^{(i)} \to X^{(i-1)}$ ($1 \leq i\leq N$) is the blow-up at a single closed point $x_{i-1}$ on $(X^{(i-1)})_s$. Let $\pi:\tX \to X$ be the composite of the above maps. Let $E_i\subset X^{(i)}$ be the exceptional divisor of $\pi_i$ and let $\wt{E_i} \subset \tX$ be its strict transform (cf.\ Appendix A below for the definition). Note that $E_i\simeq \bP^d_{x_{i-1}}$. Let $\tau_i:\wt{E_i}\to E_i$ be the natural map. We have the following facts (see Proposition \ref{blow-upncd} below for (1)--(3)):
\begin{itemize}
\item[(1)]
{\it $\wt{E_1},\dots,\wt{E_N}$ are regular irreducible components of $\tXsred$.}
\item[(2)]
{\it The reduced divisor $\bigcup_{i=1}^N \, \wt{E_i}$ on $\tX$ has simple normal crossings.}
\item[(3)]
{\it If $X\in Ob(\cQS)$, then $\tX\in Ob(\cQS)$.}
\item[(4)]
{\it For $1\leq i\leq N$, there exists a unique finite set $S_i$ of closed points on $E_i$ such that
$$
\wt{S_i}:=\tau_i^{-1}(S_i) = \wt{E_i} \cap 
\big({\bigcup}_{1\leq j \leq N,\, j \ne i} \; \wt{E_j}\big) 
$$
and such that $\wt{E_i}-\wt{S_i} \simeq E_i -S_i$ via $\tau_i$. Note that $S_i$ may be empty.}
\end{itemize}
These facts and the following lemma will be used later in \S\ref{sect8} below.
\stepcounter{thm}
\begin{lem}\label{blow-up3}
Assume $\cd_{\ell}(\rf) \le 1$. For $1\leq i\leq N$, let $W_i\subset \wt{E_i}$ be a closed integral subscheme of dimension one such that $W_i\not\subset \wt{S_i}$. Assume either that $\Lam=\Linfty$ or that $\deg_{E_i/x_{i-1}}(\tau_i(W_i))$ is prime to $\ell$ for all $1\leq i\leq N$. Then we have
$$ \Ker\big(\pi_*:\HL 2 \tX \to \HL 2 X\big) = 
\underset{1\leq i\leq N}{\sum} \Lam\langle W_i \rangle .$$
\end{lem}
\begin{pf}
Let $p$ be the projection $\tX\to X^{(1)}$. Since $\pi(W_i)$ is a closed point of $X_s$, the right hand side is contained in the left hand side by Proposition \ref{prop1-1}\,(1). There is an exact sequence
\begin{align*}
0 \lra \Ker\big(p_* : \HL 2 \tX \ra \HL 2 {X^{(1)}}\big) & \lra
\Ker\big(\pi_* : \HL 2 \tX \ra \HL 2 X \big) \\
& \os{p_*}{\lra} \Ker\big( \pi_{1*} : \HL 2 {X^{(1)}} \ra \HL 2 X \big).
\end{align*}
By induction on $N$, we may suppose that the first term is equal to $\sum_{2\leq i\leq N} \, \Lam \langle W_i \rangle$. The last group is equal to $\Lam \big\langle p_*W_1 \big\rangle$ by Corollary \ref{blow-up2}, which agrees with $p_*\big(\Lam \langle W_1 \rangle\big)$ by Proposition \ref{prop1-1}\,(1) for the proper map $p : \tX \to X^{(1)}$. Lemma \ref{blow-up3} follows easily from this fact.
\end{pf}

\bigskip

\newpage

\section{A moving lemma}\label{sect7}
Let the notation be as in \S\ref{sect1}. The discrete valuation ring $R$ is not necessarily henselian, and the residue field $\rf$ of $\ri$ is arbitrary. In this section we prove the following moving lemma for cycles on $X \in Ob(\cC)$, which will be used in the proof of Theorem \ref{injectivity}. 
\begin{prop}\label{moving}
Let $X$ be an integral regular scheme which belongs to $Ob(\cC)$. Let $Y\subset X$ be a proper closed subscheme, and put $U:=X-Y$. Then for $q \geq 0$, the natural map ${\bigoplus}_{x \in X_q \cap U}\,\bZ \to \CH_q(X)$ is surjective.
\end{prop}
\noindent
A slight modification of the following argument will show that the assertion of Proposition \ref{moving} is valid even if the base scheme $B$ is the spectrum of a Dedekind ring. Since this generalization is not used in this paper, we do not pursue this here.
\par
\bigskip
For the proof of the proposition, we need the following lemmas.

\begin{lem}\label{bertinimoving}
Let $X$ be an integral regular scheme and let $Y\subset X$ be a proper closed subscheme with $U=X-Y$. Let $y$ be a point on $Y$, and put $c:=\codim_X(y)$. Then there exists an integral closed subscheme $Z\subset X$ of codimension $c-1$ satisfying the conditions $(1)$ $Z\cap U\not=\emptyset$, $(2)$ $y\in Z$, $(3)$ $Z$ is regular at $y$. 
\end{lem}

\begin{lem}\label{bertinimoving2}
Let $Z$ be a quasi-projective scheme over $B$ and assume given a finite number of points $z_1,\dots, z_r$ on $Z$. Then there exist a dense open subset $V\subset Z$ and a dense open immersion $V\hookrightarrow U$ such that $V$ contains all $z_1,\dots, z_r$ and such that $U$ is affine.
\end{lem}

\noindent
We first prove Proposition \ref{moving} admitting these lemmas. Take an arbitrary $y\in Y_q$. Let $Z\subset X$ be as in Lemma \ref{bertinimoving}. Then $\dim(Z\cap Y)\leq q$ and $Z_q\cap Y$ is finite. Let $\pi:\tZ\to Z$ be the normalization of $Z$. The assumption that $X$ is quasi-projective over $B$ implies that $\tZ$ is quasi-projective over $B$, and hence by Lemma \ref{bertinimoving2}, there exists an open subset $V\subset \tZ$ and a dense open immersion $V\hookrightarrow U=\Spec(A)$ such that $V$ contains $\pi^{-1}(Z_q \cap Y)$. Let $A'$ be the semi-localization of $A$ at the ideals corresponding to $\pi^{-1}(Z_q \cap Y) \subset U$. It is a regular semi-local domain of dimension one, and hence a principal ideal domain. Since $Z$ is regular at $y$, there is a unique point $\tilde{y}\in \tZ$ lying over $y$. Suppose $f \in A'$ defines the principal divisor $\tilde{y}$ on $\Spec(A')$, which we regard as a rational function on $\tZ$. Then the divisor $\div_{\tZ}(f)$ on $\tZ$ has the form
$$ \div_{\tZ}(f) = [\tilde{y}] + \sum_{i=1}^{N} \ m_i \, [x_i]$$
for some points $x_i \in \tZ_q \setminus\big(\{\tilde{y}\} \cup \pi^{-1}(Z_q \cap Y)\big)$ and some integers $m_i$ for $1\leq i\leq N$. Applying $\pi_*$ to both hand sides, we get the following formula for the divisor $\div_{Z}(f)=\pi_*\div_{\tZ}(f)$ on $Z$:
$$ \div_{Z}(f) =  [y] + \sum_{i=1}^{N} \ m'_i \, [x_i']
\qquad\text{for some integers $m'_1,\dots,m'_N$},$$
where 
$x_i' =\pi(x_i)\in Z_q \setminus \big(\{y \} \cup (Z_q \cap Y)\big)$.
This proves the desired assertion.
\hfill $\square$

\bigskip\noindent
\it Proof of Lemma \ref{bertinimoving}. \rm
We have only to show the case that $Y$ is a divisor on $X$, and we may replace $X$ with $\Spec(\cO_{X,y})$. Suppose that $Y$ is defined by a non-zero element $\pi \in \cO_{X,y}$ on $X=\Spec(\cO_{X,y})$.
It suffices to show the following:
\begin{sublem}
Let $(A,\fm)$ be a regular local ring of dimension $c \ge 1$ and let $\pi\in \fm$ be a non-zero element. Then there exists a regular system of parameters $\{a_1,\dotsc,a_{c-1},a_{c}\}$ of $A$ for which $\pi$ does not belong to the ideal $(a_1,\dotsc,a_{c-1}) \subset A$.
\end{sublem}
\noindent
We prove this sublemma by induction on $c \ge 1$. The case $c=1$ is clear. Suppose $c \ge 2$.
Because $A$ is a unique factorization domain and $\fm/(\fm^2+(\pi))$ is non-zero by the assumption $c \ge 2$, one can find a prime element $a_1$ which does not divide $\pi$ and does not belong to $\fm^2$. Then $A':=A/(a_1)$ is a regular local ring of dimension $c-1$ such that the residue class $\overline{\pi}\in A'$ of $\pi$ is non-zero. By the induction hypothesis, there is a regular system of parameters $\{b_2,\dotsc,b_{c-1},b_c\}$ of $A'$ such that $\overline{\pi}\not\in (b_2,\dotsc,b_{c-1})$. Taking lifts $a_2,\dotsc,a_{c-1}, a_c\in A$ of $b_2,\dotsc,b_{c-1}, b_c$, we get a regular system of parameters $\{a_1,\dotsc,a_{c-1},a_c\}$ of $A$ such that $\pi \not\in (a_1,\dotsc,a_{c-1})$. Thus we obtain the sublemma and Lemma \ref{bertinimoving}. \hfill $\square$

\bigskip\noindent
\it Proof of Lemma \ref{bertinimoving2}. \rm
By assumption there is a locally closed immersion $i: Z\hookrightarrow \bP^N_B$, where $\bP^N_B$ denotes the projective space $\Proj(\ri [T_0,\dots,T_N])$ over $B=\Spec(\ri)$. Replacing $Z$ by the closure of $i(Z)$ in $\bP^N_B$, we may assume $Z$ is a closed subscheme of $\bP^N_B$. Then, for any affine open subscheme $U\subset \bP^N_B$, $U\cap Z$ is affine, and hence we may assume $X=\bP^N_B$. If an open subset $U\subset \bP^N_B$ contains a point in the closure of $z_i \in \bP^N_B$, then it should contain $z_i$ as well. Therefore we may assume that $z_1,\dots,z_r$ are closed points of $\bP^N_s=\bP^N_B\times_B s=\Proj(\rf[T_0,\dots,T_N])$. Let $\fm_i\subset A:=\rf[T_0,\dots,T_N]$ be the homogeneous ideal corresponding to $z_i$ for $1\leq i\leq r$. We will prove the following sublemma:
\begin{sublem}\label{bertinimoving2claim}
There exists a non-zero homogeneous polynomial $f\in A$ which is not contained in any of $\fm_1,\dotsc,\fm_r$, or equivalently, for which the hypersurface $H_f\subset \bP^N_s$ defined by $f$ does not contain any of $z_1,\dotsc,z_r$.
\end{sublem}
\noindent
We first finish the proof of Lemma \ref{bertinimoving2} admitting this sublemma. Let $f'\in \ri[T_0,\dots,T_N]$ be a homogeneous element with $f' \mod (\pi) =f$, where $f$ is taken as in the claim. Then the open subset $U=\{f'\not=0\} \subset \bP^N_B$ is affine and contains all $z_1,\dots, z_r$.

\bigskip\noindent
\it Proof of Sublemma \ref{bertinimoving2claim}. \rm
We first assume that $\rf$ is infinite.
Then the assertion follows from the following fact:
In the dual projective space $(\bP^N_s)^\vee$ (i.e. the Grassmanian of hyperplanes in $\bP^N_s$), the locus $S_i$ of those hypersurfaces that contain $z_i$ is a proper closed subset, and hence $(\bP^N_s)^\vee - (\bigcup_{1\leq i\leq r}\, S_i)$ is dense open in $(\bP^N_s)^\vee$. We next prove the case that $\rf$ is finite. By the previous case, there exists a finite Galois extension $\rf'/\rf$ and a linear form $L=\sum_{0\leq j \leq N} \, a_j T_j \in A\otimes_\rf \rf'$ such that $L\not\in \bigcup_{1\leq i \leq r}\, \fm_i\otimes_\rf \rf'$. It is easy to see that for any $\sigma\in \Gal(\rf'/\rf)$ and any $i$ with $1\leq i\leq r$, the linear form $L^\sigma=\sum_{0\leq j \leq N} \, a_j^\sigma T_j$ does not belong to $\fm_i\otimes_\rf \rf'$. Hence the homogeneous polynomial $f:=\prod_{\sigma\in \Gal(\rf'/\rf)}\, L^\sigma \; \in A$ satisfies the condition in Sublemma \ref{bertinimoving2claim}. Thus we obtain the sublemma, Lemma \ref{bertinimoving2} and Proposition \ref{moving}.
\hfill $\square$

\newpage

\section{Proof of main theorem}\label{sect8}

In this section we complete the proof of Theorem \ref{injectivity}, which consists of three steps. Assume that $R$ is excellent and henselian and that $\rf$ is finite or separably closed. We first introduce the following construction, which will be useful later:
\begin{const}\label{const8-1}
Let $X$ be an object of $\cQSP$, and suppose we are given integral closed subschemes $C_1,\dotsc,C_m \subset X$ of dimension one such that $C_j \not \subset X_s$ for any $1 \le j \le m$. Then by the resolution of singularities for embedded curves (see Theorem \ref{res1} in Appendix A below), there exists a sequence of blow-ups at single closed points as in Setting {\rm\ref{sucbu}}
\stepcounter{equation}
\begin{equation}\label{mainclaim0}
\begin{CD}
X=X^{(0)} @<{\pi_1}<< X^{(1)}  @<{\pi_2}<< X^{(2)} @<{\pi_3}<< \dotsb @<{\pi_N}<< X^{(N)} =\tX
\end{CD}
\end{equation}
such that the strict transforms $\wt{C_1},\dotsc,\wt{C_m} \subset \tX$ of $C_1,\dotsc,C_m$ are disjoint and satisfy the condition {\rm (ii)} of Theorem {\rm\ref{bertini}}. For each $1 \le i \le N$, let $E_i\subset X^{(i)}$ be the exceptional divisor of $\pi_i$ and let $\wt{E_i} \subset \tX$ be its strict transform. Let $\tau_i: \wt{E_i}\to E_i$ be the natural map. Note that $\tX$ belongs to $Ob(\cQSP)$ by Setting \ref{sucbu}\,(3).
\end{const}
\par
\smallskip
\noindent
{\bf Step 1.}
The aim of this step is to prove the following important lemma:
\stepcounter{thm}
\begin{lem}\label{mainlem1}
$\rho_X: \CH_1(X)\otimes\Linfty \to \HLinfty 2 X$
is injective for $X\in \cQSP$ with $\dim(X)=3$. 
\end{lem}
\begin{pf}
Take $\alpha\in \CH_1(X)\otimes\Linfty$ and assume $\rho_X(\alpha)=0$. By Proposition \ref{moving}, we have $\alpha = \sum_{1\leq j \leq m}\ [C_j]\otimes \lambda_j$ $(\lambda_j\in \Linfty)$ for some integral closed subschemes $C_1,\dots,C_m \subset X$ of dimension one such that $C_j\not\subset X_s$ for any $1\leq j\leq m$. Apply Construction \ref{const8-1} to $C_1,\dots,C_m$, and let $\pi:\tX \to X$ be the composite of \eqref{mainclaim0}. Applying Theorem \ref{bertini} to $\tC := \wt{C_1} \cup \dotsb \cup \wt{C_m} \hra \tX$, we obtain an ample divisor $\iota :Y \hookrightarrow \tX$ which contains $\tC$ and for which $(\tX,Y)$ is an ample \qspp. We have 
$$
\alpha=\pi_*\iota_*\beta \qwith \quad
\beta := \sum_{1\leq j \leq m} \big[\wt{C_j}\big]\otimes \lambda_j
\in \CH_1(Y)\otimes\Linfty.
$$
For each $1 \le i \le N$, let $\wt{E_i} \subset \tX$ be as in Construction \ref{const8-1}, and put $W_i:=\wt{E_i} \cap Y$. Each $W_i$ is non-empty (cf.\ Remark \ref{bertinirem}\,(3)), integral and regular of dimension one by the assumption $\dim(X)=3$, and not contained in $\wt{E_{i'}}$ for any $i' \ne i$ by Setting \ref{sucbu}\,(1) and (2). Hence by Lemma \ref{blow-up3} we have
\begin{equation}\label{H2Lblow-up0}
\Ker\big(\pi_* : \HLinfty 2 \tX \to \HLinfty 2 X \big) = \underset{1\leq i \leq N}{\sum} \Linfty \langle W_i \rangle.
\end{equation}
Now we consider the following diagram whose middle row is exact and whose square is commutative (cf.\ Propositions \ref{lem1-0}\,(2) and \ref{prop1-1}\,(1)):
\begin{equation}\label{keydiagram}
\begin{CD}
@. \CH_1(Y)\otimes\Linfty @>{\iota_*}>> \CH_1(\tX)\otimes\Linfty  \\
@. @VV{\rho_Y}V @VV{\rho_{\tX}}V \\
 \HLinfty 3 {\tX-Y} @>>> \HLinfty 2 {Y} @>{\iota_*}>> \HLinfty 2 \tX \\
\hspace{-54pt}{}_{\text{(Theorem \ref{Lcondition})}} @| @. @VV{\pi_*}V \\
 0 @. @. \HLinfty 2 X. \\
\end{CD}
\end{equation}
Let us recall the assumption that $\rho_X(\alpha)=0$. Since $\pi_*\iota_*\beta=\alpha$, we have $\pi_*\iota_*\rho_Y(\beta)=\rho_X(\pi_*\iota_*\beta)=\rho_X(\alpha)=0$ by Proposition \ref{prop1-1}\,(1) for the proper map $\pi \circ \iota : Y \to X$. Hence $\iota_*\rho_Y(\beta)$ belongs to $\Ker(\pi_*)$, and there exist $\mu_1,\cdots,\mu_N\in \Linfty$ for which the $1$-cycle
$$
\gamma:=\underset{1\leq i\leq N}{\sum}[W_i]\otimes \mu_i \in \CH_1(Y)\otimes \Linfty
$$
satisfies $\iota_*\rho_Y(\gamma)=\rho_{\tX}(\iota_*\gamma)=\iota_*\rho_Y(\beta)$, by \eqref{H2Lblow-up0} and the commutativity of the square in \eqref{keydiagram}. Then we have $\beta=\gamma$, because $\iota_*\rho_Y$ is injective by Lemma \ref{injlowdim} and Theorem \ref{Lcondition}. Since $\pi(W_i)$ is a closed point of $X$, we get $\alpha=\pi_*\iota_*\beta=\pi_*\iota_*\gamma=0$, and the proof of Lemma \ref{mainlem1} is complete.
\end{pf}

\medskip
\noindent
{\bf Step 2.} 
In this step we prove Theorem \ref{injectivity}, assuming $\dim(X)=3$.
We need the following general lemmas:
\addtocounter{thm}{3}
\begin{lem}\label{longes1}
For $X\in Ob(\cC)$ of dimension $\leq 3$ and $\Lam = \Ln$ and $\Linfty$, there is an exact sequence
$$
\HL 3 X \lra \KHL 3 X \lra \CH_1(X)\otimes \Lam \os{\rho_X}{\lra} \HL 2 X .
$$
\end{lem}
\begin{pf}
The assertion follows from the niveau spectral sequence \eqref{BOss} and \eqref{E1chow}.
\end{pf}

\begin{lem}\label{MS}
For $X\in Ob(\cC)$ of dimension $\leq 3$, we have $\KHLn 3 X = \KHLinfty 3 X [\ell^n]$.
\end{lem}
\begin{pf}
Since $\dim(X) \le 3$, we have
$$
\KHL 3 X =\Ker\left(d^1_{3,0} : {\bigoplus}_{x\in X_3} \, H^3(x,\Lam(2)) \to
{\bigoplus}_{x\in X_2} \, H^2(x,\Lam(1))\right)
$$
for $\Lam= \Ln$ and $\Linfty$, where $d^1_{3,0}$ denotes the differential map $E^1_{3,0} \to E^1_{2,0}$ of the the niveau spectral sequence \eqref{BOss}. The assertion follows from the isomorphisms
$$
H^3(x,\Ln(2))\simeq H^3(x,\Linfty(2))[\ell^n] \qaq
H^2(x,\Ln(1))\simeq H^2(x,\Linfty(1))[\ell^n],
$$
which are consequences of the divisibility of the groups
$H^2(x,\Linfty(2))$ (Merkur'ev-Suslin \cite{MS})
 and $H^1(x,\Linfty(1))$ (Hilbert's theorem 90), respectively.
\end{pf}
\noindent
We prove Theorem \ref{injectivity} assuming $\dim(X)=3$. By Theorem \ref{artin}, it is enough to prove the injectivity of $\rho_X:\CH_1(X)/\ell^n \to \HLn 2 X$. By Lemmas \ref{longes1} and \ref{MS}, the problem is reduced to showing $\KHLinfty 3 X=0$. By Theorem \ref{bertini}, we take an ample divisor $Y\subset X$ for which $(X,Y)$ is an ample \qspp. Put $U:=X-Y$. Since $\dim(X)=3$, we have $\KHLinfty 3 X \hookrightarrow \KHLinfty 3 U$ by \eqref{KCexactsequence}. We prove $\KHLinfty 3 U=0$. By Proposition \ref{prop1-1}\,(2), there is a commutative diagram with exact rows
$$
\begin{CD}
\CH_1(Y)\otimes\Linfty @>>> \CH_1(X)\otimes\Linfty @>>> \CH_1(U)\otimes\Linfty @>>> 0 \\
@VV{\rho_Y}V @VV{\rho_X}V @VV{\rho_U}V \\
\HLinfty 2 Y @>>> \HLinfty 2 X @>>> \HLinfty 2 U. \\
\end{CD}
$$
By Lemma \ref{mainlem1} and Theorem \ref{artin}, $\rho_X$ is injective and $\rho_Y$ is surjective. Hence $\rho_U$ is injective and $\CH_1(U)\otimes\Linfty$ is zero, because $\HLinfty 2 U$ is zero by Theorem \ref{Lcondition}. Finally $\HLinfty 3 U$ is zero by Theorem \ref{Lcondition}, and $\KHLinfty 3 U$ is zero as well by Lemma \ref{longes1}. This completes Step 2.
\qed
\par
\bigskip
\noindent
{\bf Step 3.}
We prove Theorem \ref{injectivity} by induction on $\dim(X)=d+1 \ge 3$. The case $d=2$ has been shown in Step 2. Assume $d \geq 3$. By Theorem \ref{artin}, it remains to prove the injectivity of $\rho_X:\CH_1(X)/\ell^n \to \HLn 2 X$.
Take $\alpha\in \CH_1(X)/\ell^n$ and assume $\rho_X(\alpha)=0$. We have to show $\alpha=0$. By Proposition \ref{moving}, we have $\alpha = \sum_{1\leq j \leq m}\ [C_j]\otimes \lambda_j$ $(\lambda_j\in \Ln)$ for some integral closed subschemes $C_1,\dots,C_m \subset X$ of dimension one such that $C_j\not\subset X_s$ for any $1\leq j\leq m$. Apply Construction \ref{const8-1} to $C_1,\dots,C_m$, and let $\tau_i : \wt{E_i} \to E_i$ ($1\leq i\leq N$) be as in Construction \ref{const8-1}. We need the following lemma, where $\rf$ is arbitrary.
\begin{lem}\label{takelines}
For each $1\leq i\leq N$, there is an integral regular closed subscheme $W_i\subset \wt{E_i}$ of dimension one which satisfies the following three conditions$:$
\begin{itemize}
\item[(1)] $\tau_i(W_i) \, (\subset E_i)$ has dimension one and degree prime to $\ell$ in the sense of Lemma {\rm\ref{blow-up3}}.
\item[(2)]
 $W_i \cap \tC$ is empty and $W_i \cap \wt{E_j}$ is empty for any $j \ne i$.
\item[(3)]
 $W_i$ meets the strict transform of $\Xsred$ $($in $\tX)$ only in its regular locus and transversally.
\end{itemize}
\end{lem}
\begin{pf}
Fix an $i$ with $1 \le i \le N$. We first prove the case that $\rf$ is infinite. Let $S_i \subset E_i$ be the finite subset as in Setting \ref{sucbu}\,(4), which may be empty. Let $C^{(i)}$ be the union of the strict transforms in $X^{(i)}$ of $C_j$'s with $j=1,\dotsc,m$. Note that $\varSigma_i:=C^{(i)} \cap E_i$ is finite, because $E_i \subset (X^{(i)})_s$ and any irreducible component of $C^{(i)}$ is not contained in $(X^{(i)})_s$. Let $Z_1,\dotsc,Z_r$ be the irreducible components of $\Xsred$ whose strict transforms $Z_1^{(i)},\dotsc,Z_r^{(i)} \subset X^{(i)}$ meet $E_i$. For $1 \le \lam \le r$, put $H_{\lam}:=Z_{\lam}^{(i)} \cap E_i$, which is a hyperplane on $E_i$. Since $\rf$ is infinite by assumption, there is a line $l_i \subset E_i$ which does not meet $S_i \cup \varSigma_i$ and intersects $H_{\lam}$ transversally for any $1 \le \lam \le r$. Then $W_i:=\tau_i^{-1}(l_i)\subset \wt{E_i}$ is the strict transform of $l_i$ in $\tX$ and we have $\tau_i|_{W_i} : W_i \isom l_i$. One can easily check that $W_i$ satisfies the conditions (1)--(3) in the lemma by the choice of $l_i \subset E_i$.

We next prove the case that $\rf$ is finite. Put $\rf_i:=\varGamma(E_i,\cO_{E_i})$, the constant field of $E_i$. Fix a prime number $\ell'$ different from $\ell$, and let $\rf_i'/\rf_i$ be the maximal pro-$\ell'$ field extension. Since $\rf'_i$ is infinite, there is a line $l_i \subset E_i \otimes_{\rf_i} \rf'_i$ which satisfies the same conditions as in the previous case with $(S_i, \varSigma_i, H_{\lam})$ replaced by $(S_i\otimes_{\rf_i} \rf'_i, \varSigma_i \otimes_{\rf_i} \rf'_i, H_{\lam}\otimes_{\rf_i} \rf'_i)$. Let $l_i' \subset E_i$ be the image of $l_i$. Then $W_i:=\tau_i^{-1}(l_i')\subset \wt{E_i}$ is the strict transform of $l_i'$ in $\tX$ and we have $\tau_i|_{W_i} : W_i \isom l_i'$. The degree of $l_i'$ in the sense of Lemma \ref{blow-up3} is a power of $\ell'$, which implies the condition (1). The conditions (2) and (3) follow from the choice of $l_i\subset E_i \otimes_{\rf_i} \rf'_i$.
\end{pf}
\par
\medskip
We turn to the proof of Theorem \ref{injectivity} and fix $W_1,\dotsc,W_N$ as in Lemma \ref{takelines}. By the condition (1) in Lemma \ref{takelines} and Lemma \ref{blow-up3} we have
\addtocounter{equation}{2}
\begin{equation}\notag
 \Ker\big(\pi_* : \HLn 2 \tX \to\HLn 2 X \big)=\underset{1\leq i\leq N}{\sum} \Ln\langle W_i \rangle.
\end{equation}
By the conditions (2) and (3) in Lemma \ref{takelines} and Theorem \ref{bertini}, there exists an ample divisor $\iota:Y \hookrightarrow \tX$ which contains $\tC$ and $\bigcup_{i=1}^N \, W_i$ and for which $(\tX,Y)$ is an ample \qspp. Here we have used the assumption $d\geq 3$ to ensure that $\dim(W_i)<\frac{1}{2}\dim\big(\wt{E_i}\big)$ (cf.\ the condition (iii) in Theorem \ref{bertini}). By the induction hypothesis, $\rho_Y:\CH_1(Y)/\ell^n \to \HLn 2 Y$ is injective. Since $d\geq 3$, we have $\HLn 3 {\tX-Y}=0$ by Theorem \ref{Lcondition}. Therefore we see that $\alpha=0$ by the same arguments as in Step 1 using the diagram \eqref{keydiagram} with $\Linfty$ replaced by $\Ln$. This completes the proof of Theorem \ref{injectivity}.
\qed
\def\index{j}

\newpage

\section{Applications of main theorem}\label{sect9}
In this section we prove Corollaries \ref{mainthmcorsepcl}\,--\,\ref{affineLintro} and Theorems \ref{mainthmfinite} and \ref{mainthmKC}. Throughout this section, $\ell$ denotes a prime number prime to $p:=\ch(\rf)$. For a scheme $X \in Ob(\cQS)$ and an effective Cartier divisor $Y$ on $X$ which belongs to $Ob(\cC)$, we define the pull-back homomorphism
$$
  \qquad  i^* : \CH_q(X) \lra \CH_{q-1}(Y) \qquad \quad (i: Y \hra X)
$$
by the same arguments as in \cite{Ful} \S2.3 and Proposition 2.6\,(a). This map sends an integral closed subscheme $f : Z \hra X$ of dimension $q$ with $Z \not \subset Y$ to the Cartier divisor $f^*Y$ on $Z$, regarded as a $(q-1)$-cycle class on $Y$. If further $X$ belongs to $Ob(\cQSP)$ and $Y$ is over $s=\Spec(\rf)$, then we define the intersection number $\langle \alpha,Y \rangle$ for $\alpha \in \CH_1(X)$ to be the value of $\alpha$ under the composite map
$$
\CH_1(X) \os{i^*}{\lra} \CH_0(Y) \os{\deg_s}{\lra} \bZ,
$$
where the last arrow is the degree map over $s$ (loc.\ cit.\ Definition 1.4). Note that the degree map is well-defined by the properness of $Y$. We first show
\begin{prop}\label{lemcor0}
Let $X$ be an object of $\cQSP$, and let $C \subset X$ be an integral closed subscheme of dimension one which is flat over $B$. Let $x$ be the closed point of $C$. Let $Y$ be an irreducible component of $\Xsred$ which contains $x$, and let $i$ be the closed immersion $Y \hra X$. Put $d:=\dim(X)-1$. Let $\cl_X(-) \in H^{2d}(X,\Ln(d))$ and $\cl_{Y}(-) \in H^{2d}(Y,\Ln(d))$ be the cycle class of $1$-cycles on $X$ and $0$-cycles on $Y$, respectively {\rm(}\cite{Fu} Definition {\rm1.1.2}, see also loc.\ cit.\ Corollary {\rm1.1.5)}. Then{\rm:}
\begin{enumerate}
\item[(1)] We have $i^*\cl_X(C)=\cl_Y(i^*C)$ in $H^{2d}(Y,\Ln(d))$.
\item[(2)] $\langle C,Y \rangle$ is divisible by $\deg_s(x)$ and we have
\begin{equation}\notag
  i^* \cl_X(C) = \frac{\langle C,Y \rangle}{\deg_s(x)} \cdot \cl_{Y}(x) \;\; \hbox{ in } \;\; H^{2d}(Y,\Ln(d)).
\end{equation}
\end{enumerate}
\end{prop}
\begin{pf}
Put $\vare:=C \times_X Y$, which is the spectrum of an artinian local ring $A$ and represents the $0$-cycle $i^*C$ on $Y$. We have $\langle C,Y \rangle={\mathrm{length}}_A(A) \cdot \deg_s(x)$ by the definition of the intersection number $\langle -,Y \rangle$. We have $i^* \cl_X(C)= \cl_Y(\vare)$ by \cite{Fu} Proposition 1.1.9, and $\cl_Y(\vare) ={\mathrm{length}}_A(A) \cdot \cl_Y(x)$ by \cite{SGA4.5} Cycle Th\'eor\`eme 2.3.8\,(i). The assertions follow from these facts.
\end{pf}
\par
\medskip
We prove Corollaries \ref{mainthmcorsepcl}\,--\,\ref{affineLintro}.
\begin{cor}[{{\bf Corollary \ref{mainthmcorsepcl}}}]\label{mainthmcorsepcls9}
Let $X$ be an object of $\cQSP$, and let $Y_1,\dots,Y_N$ be the irreducible components of $\Xsred$. Assume that $\rf$ is separably closed. Then the assignment $\alpha \in \CH_1(X) \mapsto \langle \alpha,Y_{\index} \rangle_{1\leq \index \leq N} \in \bZ^N$ induces an isomorphism
$$
 \CH_1(X)/\ell^n \os{\simeq}{\lra} \bigoplus_{1\leq \index \leq N}\, \bZ/\ell^n \bZ
   \quad \hbox{ for any } \, n \ge 1.
$$
\end{cor}
\begin{pf}
Put $d:=\dim(X)-1$. We note the following isomorphisms analogous to \eqref{H1isom}:
\addtocounter{equation}{2}
\stepcounter{thm}
\begin{equation}\label{corisom1}
\HLn 2 X \os{\eqref{purity}}{\simeq} \Het {2d}(X,\Ln(d)) \isom
\bigoplus_{1\leq \index \leq N}\, \Het {2d}(Y_{\index},\Ln(d)) \isom \underset{1\leq \index \leq N}{\bigoplus}\Ln,
\end{equation}
where the second arrow is bijective by a similar argument as for Lemma \ref{sublemH1isom}\,(2). The last arrow is the direct sum of the trace maps for $Y_{\index} \to s$ (\cite{SGA4} XVIII Th\'eor\`eme 2.9), which is bijective by \cite{Mi} VI Lemma 11.3. Hence the corollary follows from Theorem \ref{injectivity}, Proposition \ref{moving} and the following lemma:
\begin{lem}\label{lemcor1}
Let $C \subset X$ be an integral closed subscheme of dimension one which is flat over $B$. Then the composite map \eqref{corisom1} sends the cycle class $\rho_X(C)$ to $\langle C,Y_{\index} \rangle_{1\leq \index \leq N}$.
\end{lem}
\medskip
\noindent
{\it Proof of Lemma \ref{lemcor1}.}
Let $x$ be the closed point of $C$. By Lemma \ref{lem1-1}\,(1), the first isomorphism in \eqref{corisom1} sends $\rho_X(C)$ to the cycle class $\cl_X(C)$. By Proposition \ref{lemcor0}\,(2), we have
\begin{equation}\notag
  i_{\index}^* \cl_X(C) = \frac{\langle C,Y_{\index} \rangle}{\deg_s(x)} \cdot \cl_{Y_{\index}}(x) \;\; \hbox{ in } \;\; \Het {2d}(Y_{\index},\Ln(d))
\end{equation}
for $1 \le \index \le N$. Here $i_{\index}$ denotes the closed immersion $Y_{\index} \hra X$, and $\cl_{Y_{\index}}(x)$ means zero if $x \not\in Y_{\index}$. The assertion of the lemma follows from this formula. Indeed if $x \in Y_{\index}$, then the trace map $\Het {2d}(Y_{\index},\Ln(d)) \to \Ln$ sends $\cl_{Y_{\index}}(x)$ to $\deg_s(x)$ by \cite{SGA4.5} Cycle Th\'eor\`eme 2.3.8\,(i). This completes the proof of Lemma \ref{lemcor1} and Corollary \ref{mainthmcorsepcls9}.
\end{pf}
\begin{cor}[{{\bf Corollary \ref{goodreduction}}}]\label{goodreductions9}
Assume that $X$ is smooth and projective over $B$ and that $\rf$ is finite or separably closed. Then for any $n \ge 1$, we have
$$s
\begin{CD}
\CH_0(X_s)/\ell^n @<{\simeq}<{i^*}< \CH_1(X)/\ell^n @>{\simeq}>{j^*}> \CH_0(X_\eta)/\ell^n \quad (X_s \, \os{i}{\hra} \, X \, \os{j}{\hookleftarrow} \, X_\eta).
\end{CD}
$$
Here $\CH_0(X_\eta)$ denotes the Chow group of $0$-cycles on $X_{\eta}$ in the usual sense {\rm(}cf.\ Remark {\rm\ref{rem1-1})}.
\end{cor}
\begin{pf}
By Lemma \ref{bertinimoving} and Proposition \ref{lemcor0}\,(1), there is a commutative diagram
$$
\begin{CD}
\CH_1(X)/\ell^n @>{\cl_X}>> \Het {2d}(X,\Ln(d)) \\ @V{i^*}VV @VV{i^*}V \\
\CH_0(X_s)/\ell^n @>{\cl_{X_s}}>> \Het {2d}(X_s,\Ln(d)). \\
\end{CD}
$$
where $\cl_X$ and $\cl_{X_s}$ send cycles to their cycle classes. See \cite{JSS} Theorem 1.1.1 for why these assignments factor through the Chow groups. The right vertical $i^*$ is bijective by \cite{SGA4} XII Corollaire 5.5\,(iii), and $\cl_X$ is bijective by Theorem \ref{injectivity} and Lemma \ref{lem1-1}. The arrow $\cl_{X_s}$ is obviously bijective if $\rf$ is separably closed. If $\rf$ is finite, then $\cl_{X_s}$ is bijective by the unramified class field theory for the smooth and projective $F$-variety $X_s$ (\cite{KS}, \cite{CTSS} p.\ 792 Th\'eor\`eme 5, see also loc.\ cit.\ p.793 Remarque 3). Hence the left vertical $i^*$ is bijective.
\par
 To show the bijectivity of $j^*$ in the corollary, we note the following exact sequence, which follows from the same argument as for \cite{Ful} Proposition 1.8:
$$
\begin{CD}
\CH_1(X_s)/\ell^n \os{i_*}{\lra} \CH_1(X)/\ell^n \os{j^*}{\lra} \CH_0(X_\eta)/\ell^n \lra 0.
\end{CD}
$$
We show the map $i_*$ is zero. Because $X_s$ is a principal divisor on $X$, the composite map
$$
\CH_1(X_s)/\ell^n \os{i_*}{\lra} \CH_1(X)/\ell^n \us{\simeq}{\os{i^*}{\lra}} \CH_0(X_s)/\ell^n
$$
is zero by \cite{Ful} Proposition 2.6\,(c). Hence the bijectivity of $i^*$ implies that $i_*$ is zero. Thus we obtain Corollary \ref{goodreductions9}.
\end{pf}
\begin{cor}[{{\bf Corollary \ref{affineLintro}}}]\label{affineLs9}
Let $(X,Y)$ be an ample \qspp~ with $\dim(X)=d+1\ge 2$, and let $i:Y \hookrightarrow X$ be the natural closed immersion. Assume that $\rf$ is finite {\rm(}resp.\ separably closed{\rm)}. Then for any $n \ge 1$, the push-forward map $i_* : \CH_1(Y)/\ell^n \to \CH_1(X)/\ell^n$ is bijective for $d \ge 3$ and surjective for $d=2$ {\rm(}resp.\ bijective for $d \ge 2$ and surjective for $d=1${\rm)}.
\end{cor}
\begin{pf}
Put $\Lam := \Ln$. By Propositions \ref{lem1-0}\,(2) and \ref{prop1-1}\,(2), there is a commutative diagram with exact lower row
\[
\begin{CD}
@. \CH_1(Y)/\ell^n @>{i_*}>> \CH_1(X)/\ell^n  \\
@. @V{\rho_Y}V{\simeq}V @V{\rho_X}V{\simeq}V \\
\HLn 3 U @>>> \HLn 2 Y @>{i_*}>> \HLn 2 X @>>> \HLn 2 U, \\
\end{CD}
\]
where the vertical arrows are bijective by Theorem \ref{injectivity}. Corollary \ref{affineLintro} follows from this diagram and the vanishing results in Theorem \ref{Lcondition}.
\end{pf}

\medskip
Next we prove the following result:
\begin{thm}[{{\bf Theorem \ref{mainthmfinite}}}]\label{mainthmfinites9}
Assume that $\ri$ is a $p$-adic integer ring. Let $X$ be an object of $\cQSP$, and put $V:= X_{\eta}$. Then$:$
\begin{enumerate}
\item[$(1)$]
The group $ A_0(V) \hatotzl :=\varprojlim_{n \in \mathbb N} \ A_0(V)/\ell^n$ is finite for any prime $\ell\not=p$ and is zero for almost all primes $\ell$.
\item[$(2)$]
Conjecture {\rm\ref{mainconjfinite}} holds true for $V$.
\end{enumerate}
\end{thm}
\begin{pf}
We first show (1). By Theorem \ref{bertini0} and Corollary \ref{affineLs9}, the problem is reduced to the case $\dim(V)=1$. Let $J$ be the Jacobian variety of $V$. Then $A_0(V)$ is a subgroup of $J(\qf)$ of finite index, and the assertion follows from a theorem of Mattuck \cite{Ma} and \cite{JS1} Lemma 7.8.
\par
The assertion (2) follows from (1) and loc.\ cit.\ Lemma 7.7 applied to
\begin{align*}
 &(A,\, \{B_n\}_{n \in \bN, (n,p)=1}, \, \{\varphi_n:A/n \to B_n\}_{n \in \bN, (n,p)=1}) \\
 & = (A_0(V), \, \{A_0(V)/n\}_{n \in \bN, (n,p)=1}, \, \{\id : A_0(V)/n \to A_0(V)/n\}_{n \in \bN, (n,p)=1}).
\end{align*}
This completes the proof of Theorem \ref{mainthmfinites9}.
\end{pf}

\medskip
Finally we prove Theorem \ref{mainthmKC}. The case $a\leq 2$ has been shown in Proposition \ref{mainthmKClowdim} and Theorem \ref{artincor}. Let $\Lam$ denote $\Ln$ if $\rf$ is separably closed, and denote $\Linfty$ if $\rf$ is finite. Let $X$ be an object of $\cQSP$. We have to show $\KHL 3 X=0$. We proceed to the proof by induction on $\dim(X)$. The assertion is clear if $\dim(X)\leq 2$. Assume $\dim(X) \ge 3$. By Theorem \ref{bertini} we take an ample divisor $Y\subset X$ such that $(X,Y)$ is a \qspp, and put $U:=X-Y$. By the exact sequence of complexes in \eqref{KCexactsequence}, there is a long exact sequence
\[ \dotsb \lra \KHL 3 Y \lra \KHL 3 X \lra \KHL 3 U \lra \dotsb. \]
By the induction hypothesis, it is enough to show $\KHL 3 U$ is zero. By \eqref{E1chow} and Lemma \ref{E1vanish}, the niveau spectral sequence \eqref{BOss} yields an exact sequence
\[ \HL 3 U \lra \KHL 3 U \lra \CH_1(U)\otimes\Lam. \]
By Theorem \ref{Lcondition} we have $\HL 3 U=0$. By Corollary \ref{affineLs9} and Proposition \ref{prop1-1}\,(2), we have $\CH_1(U)\otimes\Lam=0$. Hence $\KHL 3 U$ is zero. This completes the proof of Theorem \ref{mainthmKC}.
\hfill \qed

\appendix

\def\Bl{{\mathrm{Bl}}}
\def\gr{{\mathrm{gr}}}

\newpage

\section{Resolution of singularities for embedded curves}
\begin{center}
by {\sc uwe jannsen}
\par
\bigskip
{\small
Fakult\"at f\"ur Mathematik,
Universit\"at Regensburg
}
\end{center}
\par
\vspace{-2mm}
\begin{center}
{\scriptsize
Universit\"atsstr.\ 31, 93040 Regensburg, GERMANY
}
\end{center}
\vspace{-2mm}
\begin{center}
{\scriptsize
uwe.jannsen@@mathematik.uni-regensburg.de
}
\end{center}
\par
\bigskip
The purpose of this appendix is to show the following form of resolution of singularities for curves embedded in a regular scheme, which has been used in the proof of Theorem \ref{injectivity} of the main body (cf.\ \S\ref{sect8}). We believe that it has been known to experts but there is no written proof of it as far as we know.
\par
\bigskip
For a noetherian scheme $X$ and a closed subscheme $Z \subset X$, let
$$
\pi : \Bl_Z(X)=\Proj\left({\bigoplus}_{n\geq 0} \, \cI_Z^n \right) \lra X
$$ 
denote the blow-up of $X$ in $Z$, i.e., the blow-up of the ideal sheaf $\cI_Z \subset \cO_X$ defining $Z$, and let $E_Z(X) = Z\times_X\Bl_Z(X)$ denote the exceptional divisor. If $i: Y \hookrightarrow X$ is another closed immersion, then the strict transform of $Y$ in $\Bl_Z(X)$ is by definition $\Bl_{Z\times_X Y}(Y)$, the blow-up of the inverse image ideal sheaf $i^{-1}\cI_Z \cdot \cO_Y \subset \cO_Y$. By the universal property of blow-ups there is a canonical morphism $\Bl_{Z\times_X Y}(Y) \to \Bl_Z(X)$, which is a closed immersion.

\begin{thm}\label{res1}
Let $X$ be a regular excellent scheme and let $Y$ be a reduced divisor with simple normal crossings on $X$ {\rm(}cf.\ Definition {\rm\ref{defqss}\,(1))}. Let $W$ be a one-dimensional reduced closed subscheme of $X$ such that any irreducible component of $W$ is not contained in $Y$. Then there exists a proper morphism $\pi: X' \rightarrow X$ which satisfies the following conditions$:$
\begin{itemize}
\item[{\rm (i)}]
$\pi$ is the composite of a sequence of blow-ups $X'=X^{(n)} \to \dotsb \to X^{(1)} \to X$ at closed points on $W$ and closed points on the strict transforms of $W$ in $X^{(i)}$ $(1 \le i \le n-1)$.
\item[{\rm (ii)}]
The closed subset $\pi^{-1}(Y)$ with reduced subscheme structure is a reduced divisor with simple normal crossings.
\item[{\rm (iii)}]
The strict transform $W' \subset X'$ of $W$ is regular and intersects $Y'$ transversally, which means that $W'$ intersects $Y'$ only in the regular locus of $Y'$ and transversally.
\end{itemize}
\end{thm}
\par
\bigskip
By Proposition \ref{blow-upncd}\,(3) below, the condition (ii) in Theorem \ref{res1} is a consequence of the condition (i). The proof of the theorem is achieved in several steps. The following statement would be well-known, but we include a proof for lack of a suitable reference.

\begin{prop}\label{resbl-1}
Let $X$ be a reduced excellent scheme of dimension one. Then after a sequence of blow-ups $X^{(n)}\to \cdots \to X^{(1)} \to X$ at singular closed points, the resulting scheme $X^{(n)}$ is regular.
\end{prop}

\begin{pf}
Since $X$ is reduced and excellent, it has a regular dense open subset by \cite{EGAIV} 7.8.6\,(iii). Because blow-ups are compatible with Zariski localization, we may assume that $X$ is an affine scheme $\Spec(R)$. Let $x$ be a singular closed point of $X$, and let $\fm \subset R$ be the maximal ideal that corresponds to $x$. Let $\pi$ be the blow-up of $X$ at $x$:
$$
\pi: X' :=\Bl_x(X)= \Proj(S) \longrightarrow X \;\; \mbox{ with } \;\; S := \underset{n\geq 0}{\bigoplus} \; \fm^n.
$$
We first show
\begin{lem}\label{resbl-1-1}
$\pi$ is a finite morphism. In particular, $X'$ is affine by \cite{EGAII} {\rm 6.7.1}.
\end{lem}
\begin{pf*}{Proof of Lemma \ref{resbl-1-1}}
Since $\pi$ is projective, it is enough to show that $\pi$ is quasi-finite (cf.\ \cite{Mi} I Theorem 1.8). Since $\pi$ is an isomorphism over $X\setminus \{x\}$, it suffices to show that the fiber
$$
\pi^{-1}(x) = \Proj(\ol S) \;\; \hbox{ with } \;\; \ol S :=S \otimes_R R/\fm = \underset{n\geq 0}{\bigoplus}\; \fm^n/\fm^{n+1}
$$
consists of finite number of points. Let $R_{\fm}$ be the local ring $\cO_{X,x}$, and let $\fm R_{\fm}$ be its maximal ideal. Put $k := R/\fm=R_{\fm}/\fm R_{\fm}$. There is an isomorphism of finitely generated $k$-algebras
$$
\overline S \, \simeq \, \underset{n\geq 0}{\bigoplus}\; (\fm R_{\fm})^n/(\fm R_{\fm})^{n+1} =: \gr^{\bullet} R_{\fm}.
$$
Because $R_{\fm}$ is one-dimensional by assumption, $\gr^{\bullet} R_{\fm}$ has dimension one as well by theory of noetherian local rings. Hence $\pi^{-1}(x) = \Proj(\gr^{\bullet} R_{\fm})$ has dimension zero by \cite{Ku} II Proposition 4.4\,(b), and the lemma follows.
\end{pf*}
\par
\medskip
We turn to the proof of Proposition \ref{resbl-1}. Let $X_1$, \dots, $X_N$ be the irreducible components of $X=\Spec(R)$. Since $X$ is reduced, we may write $X_i=\Spec(R_i)$ $(1\leq i\leq N$), where $R_i$ is an integral domain of dimension one. Let $\wt{R_i}$ be the normalization of $R_i$ in its fraction field $K_i$. Put $K:=K_1\times \dotsb \times K_N$ and $\tR := \wt{R_1} \times \dotsb \times  \wt{R_N}$. We want to show that after a sequence of blow-ups $X^{(n)}\to \cdots \to X^{(1)}=X'\to X$ at singular points, we get $X^{(n)}=\Spec(\tR)$. By Lemma \ref{resbl-1-1} we have $X' = \Spec(R_1)$ for some finite $R$-algebra $R_1$.  Since $X$ is reduced by assumption, $R_1$ is reduced. Hence we have inclusions $R \subset R_1 \subset \tR \subset K$. As $R$ is excellent, $\tR$ is finite over $R$ by \cite{EGAIV} 7.8.6\,(ii). As $\dim(R)=1$, $\tR/R$ has finite length over $R$ (it is supported in the finitely many singular points of $\Spec(R)$). By \cite{EGAII} 8.1.7, the ideal $\fm R_1\subset R_1$ is invertible, while $\fm \subset R$ is not invertible because $x$ is a singular point of $\Spec(R)$ and $\fm R_{\fm}\subset R_{\fm}$ is not principal. Therefore $R_1\neq R$, and $\ell_{R_1}(\tR/R_1) \leq \ell_R(\tR/R_1) < \ell_R(\tR/R)$, where $\ell_R(-)$ denotes the length over $R$. Hence we obtain Proposition \ref{resbl-1} by induction on $\ell_R(\tR/R)$.
\end{pf}

\begin{prop}\label{blow-upncd}
Let $X$ be an integral regular noetherian scheme of dimension $d$, and let $Z \subset X$ be an integral regular closed subscheme with $Z\neq X$. Then$:$
\begin{itemize}
\item[(1)]
$\Bl_Z(X)$ is again integral and regular of dimension $d$, and the exceptional divisor $E_Z(X)$ is integral and regular.
\item[(2)]
Let $Y_1, Y_2 \subset X$ be regular closed subschemes containg $Z$. Assume that $Y_1$ and $Y_2$ intersect transversally. Then the canonical morphism
$$
\Bl_Z(Y_1\times_XY_2) \lra \Bl_Z(Y_1)\times_{\Bl_Z(X)}\Bl_Z(Y_2)
$$
is an isomorphism.
\item[(3)]
Let $Y \subset X$ be a reduced divisor with simple normal crossings on $X$ and let $Y_1, \dotsc, Y_r$ be the irreducible components of $Y$. Assume that all $Y_i$ contain $Z$. Then the union of $E_Z(X)$ and the strict transforms of $Y_i$ in $\Bl_Z(X)$ for $1 \le i \le r$ is a reduced divisor with simple normal crossings on $\Bl_Z(X)$.
\end{itemize}
\end{prop}

\begin{pf}
We first show (1). Since $X$ is integral and $Z \ne X$, $\Bl_Z(X)$ is integral of dimension $d$ by \cite{EGAII} 8.1.4. Since $X$ and $Z$ are regular, $\Bl_Z(X)$ is regular by \cite{EGAIV} 19.4.3 and 19.4.4. It follows from loc.\ cit.\ 19.4.2, that $E_Z(X)$ is isomorphic to the projective bundle over $Z$ associated to the conormal sheaf of $Z$ in $X$. Hence $E_Z(X)$ is integral and regular. To prove Proposition \ref{blow-upncd}\,(2), we need the following standard fact (cf.\ \cite{Ful} Remark A.6).
\begin{lem}\label{explicitblow-up}
Let $A$ be a commutative ring with unity, and let $a_1,\dotsc,a_m \in A$ be a regular sequence. Put $I :=(a_1,\dotsc,a_m) \subset A$. Then there is an isomorphism of graded $A$-algebras
$$
A[T_1,\dots,T_m]/(a_iT_j-a_jT_i \;\vert\;  1\leq i<j \leq m)\,\isom
\underset{n\geq 0}{\bigoplus}\, I^n, \;\; T_i \mapsto a_i \in I.
$$
\end{lem}

We start the proof of Proposition \ref{blow-upncd}\,(2). Since the problems are local in the Zariski topology, we may assume that $X = \Spec(A)$ for a noetherian regular local ring $A$ with maximal ideal $\fm$. Put $d := \dim(A)$. By the assumptions on $Y_1,Y_2$ and $Z$, there is a regular system of parameters \{$a_1,\dotsc,a_d\}$ of $A$ such that $Z=\Spec(A/I)$ and $Y_i=\Spec(A/I_i)$ ($i=1,2$), where
$$
I = (a_1,\dotsc,a_m),\quad
I_1=(a_1,\dotsc,a_\mu),\quad
I_2=(a_{\mu+1},\dotsc,a_{\mu+\nu})
$$
with $\mu+\nu\leq m:= \codim_X(Z)$. By Lemma \ref{explicitblow-up}, we have
$$
\Bl_Z(X) \simeq \Proj(B),\quad
\Bl_Z(Y_i) \simeq \Proj(B_i) \;\;(i=1,2),\quad
\Bl_Z(Y_1\times Y_2) \simeq \Proj(B_{12}),
$$
where $B := A[T_1,\dots,T_m]/( a_iT_j-a_jT_i \;\vert\; 1\leq i<j \leq m)$ and
{\allowdisplaybreaks
\begin{align*}
B_1 & := (A/I_1)[T_{\mu+1},\dotsc,T_m]/( a_iT_j-a_jT_i \,\vert\, \mu+1 \leq i<j \leq m) \\ & = B/( T_1,\dotsc,T_\mu,a_1,\dotsc,a_\mu),\\
B_2 & := (A/I_2)[T_1,\dotsc,T_\mu,T_{\mu+\nu+1},\dotsc,T_m]/(a_iT_j-a_jT_i \,\vert\, i,j\in [1,\mu] \cup [\mu+\nu+1, m]) \\ & = B/( T_{\mu+1},\dotsc,T_{\mu+\nu},a_{\mu+1}.\dotsc,a_{\mu+\nu}), \\
B_{12} & := \big(A/(I_1+I_2)\big)[T_{\mu+\nu+1},\dotsc,T_m]/(a_iT_j-a_jT_i\, \vert\, \mu+\nu+1 \leq i<j \leq m) \\ & = B/( T_1,\dotsc,T_{\mu+\nu},a_1,\dotsc,a_{\mu+\nu}).
\end{align*}
}
\hspace{-7pt}
The isomorphism in question is a direct consequence of these descriptions.
\par
Finally we prove Proposition \ref{blow-upncd}\,(3). By (1), the exceptional divisor $E_Z(X) \subset \Bl_Z(X)$ and the strict transforms $\Bl_{Z}(Y_i)$ of $Y_i$ ($1\le i \le r$) are regular divisors on $\Bl_Z(X)$. For $p\geq 1$ and $1\leq i_1< i_2 < \dotsb < i_p \leq s$, let
$$
Y_{i_1,\dotsc,i_p} = Y_{i_1}\times_X \dotsb \times_XY_{i_p}
$$
be the $p$-fold intersection of irreducible components of $Y$.
Applying (2) repeatedly, we have
$$
 \Bl_Z(Y_{i_1}) \times_{\Bl_Z(X)}\dotsb \times_{\Bl_Z(X)} \Bl_Z(Y_{i_p}) \simeq \Bl_Z(Y_{i_1,\dotsc,i_p}),
$$
where the scheme on the right hand side is regular by (1). This shows that the reduced divisor $\bigcup_{i=1}^r \Bl_Z(Y_i)$ has simple normal crossings. As for the intersection of $\Bl_Z(Y_{i_1,\dotsc,i_p})$ and $E_Z(X) = Z\times_X\Bl_Z(X)$, we have
\begin{align*}
E_Z(X)\times_{\Bl_Z(X)}\Bl_Z(Y_{i_1,\dotsc,i_p}) =& Z \times_X \Bl_Z(X)\times_{\Bl_Z(X)}\Bl_Z(Y_{i_1,\dotsc,i_p}) \\
=& Z\times_X \Bl_Z(Y_{i_1,\dotsc,i_p}) = Z \times_{Y_{i_1,\dotsc,i_p}} \Bl_Z(Y_{i_1,\dotsc,i_p}),
\end{align*}
which is the exceptional divisor of the blow-up $\Bl_Z(Y_{i_1,\dotsc,i_p})\rightarrow Y_{i_1,\dotsc,i_p}$ and is a regular divisor on $\Bl_Z(Y_{i_1,\dotsc,i_p})$ by (1). This shows that $E_Z(X)$ intersects all $Y_{i_1,\dotsc,i_p}$'s transversally. This completes the proof of Proposition \ref{blow-upncd}.
\end{pf}

\begin{prop}\label{curveanddivisor}
Let $X$ be a regular noetherian scheme, let $D$ be a union of finitely many regular divisors on $X$, and let $C \subset X$ be an integral regular closed subscheme of dimension one. Assume that $C$ is not contained in $D$. Let $X^{(1)}$ be the blow-up of $X$ at the $($closed\,$)$ points of $C \cap D$. Let $C^{(1)}$ be the strict transform of $C$ in $X^{(1)}$, and let $D^{(1)}$ be the union of the strict transform of $D$ and the exceptional divisors of the blow-up $X^{(1)} \to X$. For $n\geq 2$, define $(X^{(n)},C^{(n)},D^{(n)})$ inductively by applying this procedure to $(X^{(n-1)},C^{(n-1)},D^{(n-1)})$. Then $C^{(n)}$ intersects $D^{(n)}$ only in the regular locus of $D^{(n)}$ and transversally, for $n$ large enough.
\end{prop}

\begin{pf}
Because $C\cap D$ is discrete and the question is local around the points $x\in C\cap D$, we may localize $X$ to assume $X=\Spec(A)$ for a noetherian regular local ring $A$ and that $x$ corresponds to the maximal ideal $\fm$ of $A$. Suppose that $A$ has dimension $d$, and take a regular system of parameters $\{a_1,\dotsc, a_d\}$ of $A$ for which $C=\Spec(\overline A)$ with ${\ol A} := A/(a_1,\dotsc,a_{d-1})$. Let $D_1,\dots, D_s$ be the irreducible components of $D$. For $1\leq i\leq s$, we have $D_i=\Spec\big(A/(f_i)\big)$ for some $f_i\in A$. The assumption that $C$ is not contained in $D$ implies that $f_i\notin (a_1,\dotsc,a_{d-1})$ for any $i$. We now fix $i\in [1,s]$. By assumption, ${\ol A} =A/(a_1,\dotsc,a_{d-1})$ is a discrete valuation ring with a prime element $a_d \mod (a_1,\dotsc,a_{d-1})\in {\ol A}$. We define an integer $r_x(C,D_i)$ as the (normalized) valuation of $f_i \mod (a_1,\dotsc,a_{d-1})\in {\ol A}$. It depends only on the triple $(C,D_i,x)$ and not on the choice of a regular system of parameters. We will prove
\begin{lem}\label{curveanddivisorclaim}
Let $\pi: \tX=\Bl_x(X) \to X $ be the blow-up of $X=\Spec(A)$ at $x$, and let $E \subset \tX$ be the exceptional divisor. Let $\tC \subset \tX$ $($resp.\ $\wt{D_i}\subset \tX)$ be the strict transform of $C$ $($resp.\ $D_i)$. Then$:$
\begin{itemize}
\item[(1)]
$\tC$ and $E$ meet transversally.
\item[(2)]
If $r_x(C,D_i)=1$, then $\tC \cap \wt{D_i}$ is empty.
\item[(3)] 
If $r_x(C,D_i)\geq 2$, then $r_y(\tC,\wt{D_i}) < r_x(C,D_i)$ for any point $y \in \tC \cap \wt{D_i}$.
\end{itemize}
\end{lem}
\par
\smallskip
\noindent
These statements imply that for a sufficiently large $n$, the strict transform $C^{(n)}$ of $C$ in $X^{(n)}$ meets only the exceptional divisor $E^{(n)}$ of the last blow-up among the irreducible components of the inverse image of $D$. Moreover, $C^{(n)}$ intersects $E^{(n)}$ transversally. This implies Proposition \ref{curveanddivisor}. It remains to show the lemma.
\par
\bigskip
\noindent
{\it Proof of Lemma \ref{curveanddivisorclaim}.}
We show (1). By Lemma \ref{explicitblow-up} we have
\begin{align*}
\tX & = \Proj(B) \;\;\hbox{ with } \;\;
    B=A[T_1,\dotsc,T_d]/(a_iT_j-a_jT_i~\vert~1\leq i<j\leq d)\\
\tC & = \Proj(B') \;\;\hbox{ with } \;\; B'= \ol A \,[T_d],
\end{align*}
where $\tC$ is of course isomorphic to $C$. Moreover, the closed immersion $\tC \hra \tX$ is induced by the ring homomorphism $B \to B'$ sending $T_i \mapsto 0$ for $i \ne d$. This shows that $\tC$ lies in the chart $V_+(T_d)=\{T_d \neq 0\}\subset \Proj(B)$, which is the spectrum of the ring
\addtocounter{equation}{7}
\stepcounter{thm}
\begin{equation}\label{chartUd}
A[t_1,\dotsc,t_{d-1}]/(a_d\,t_i-a_i~\vert~i=1,\dotsc,d-1)\,\quad (t_i = T_i/T_d)
\end{equation}
In this chart, $E\cap V_+(T_d)$ (resp.\ $\tC\cap V_+(T_d)$) is described as the locus where  $a_d=0$ (resp. $t_1 =\dotsb = t_{d-1} = 0$). This immediately implies Lemma \ref{curveanddivisorclaim}\,(1).

To show Lemma \ref{curveanddivisorclaim}\,(2) and (3), we need to compute the strict transform $\wt{D_i}$ of $D_i$. For simplicity we put $D:=D_i$, $f:=f_i$ and $r:=r_x(C,D)$. We normalize the regular system $\{a_1,\dotsc, a_d\}$ as follows. We have $f = u \cdot a_d^r + g$ for some $g \in (a_1,\dotsc,a_{d-1})$ and some $u \in A^{\times}$. Replacing $f$ by $u^{-1}f$, we suppose $f = a_d^r + g$. We distinguish two cases. If $r=1$, then $D$ intersects $C$ transversally. Replacing $a_d$ with $f$, we suppose
$$
f = a_d.
$$
If $r\geq 2$, then we have $f = a_d^r + a_1g_1 + \dotsb + a_{d-1}g_{d-1}$ for some $g_j \in A$. Note that not all $g_j$'s are contained in $\fm$ since $f \notin \fm^2$ by regularity of $D$. If $g_k \notin \fm$, then $g_k$ is a unit and we replace $a_k$ with $a_1g_1 + \dotsb + a_{d-1}g_{d-1}$. By renumbering if necessary, we suppose
$$
f = a_d^r + a_1.
$$
Because $\tC$ is contained in the affine chart $V_+(T_d) \subset \tX$, we compute $\wt D \cap V_+(T_d)$. Note that $\wt D$ is an integral closed subscheme of $\pi^{-1}(D)=D\times_X \tX$, and that $\pi^{-1}(D) \cap V_+(T_d)$ is defined by $f$, regarded as an element of the ring \eqref{chartUd}. When $r=1$, we assumed $f=a_d$, which defines $E$ on $V_+(T_d)$. Since $\wt{D} \not\subset E$, this implies $\wt D \cap V_+(T_d)=\emptyset$ and hence $\wt D$ does not meet $\tC$. This shows Lemma \ref{curveanddivisorclaim}\,(2). Now consider the case $r\geq 2$, where we assumed $f = a_d^r + a_1$. In the ring \eqref{chartUd}, we have a factorization
$$
 f = a_d(a_d^{r-1} + t_1),
$$
which implies that the divisor $\wt{D} \cap V_+(T_d)\subset V_+(T_d)$ is defined by $a_d^{r-1} + t_1$. Since $\tC$ is defined by the equations $t_1=\dotsb=t_{d-1}=0$, we see that $\wt{D}$ meets $\tC$ only at the point $y \in V_+(T_d)$ corresponding to the ideal $(a_d,t_1,\dotsc,t_{d-1})$, and that $r_y(\tC,\wt D)=r-1$. This shows Lemma \ref{curveanddivisorclaim}\,(3), and the proof of Proposition \ref{curveanddivisor} is complete.
\end{pf}

\medskip
We deduce Theorem \ref{res1} from Propositions \ref{resbl-1}, \ref{blow-upncd} and \ref{curveanddivisor}. Let $Y\subset X$ and $W\subset X$ be as in the theorem. Since $X$ is excellent by assumption, $Z$ is excellent as well by \cite{EGAIV} 7.8.6\,(i). After a sequence of blow-ups $X^{(n)} \to \dotsb \to X^{(1)} \to X$ at singular points of $W$ and singular points of the strict transform $W^{(i)} \subset X^{(i)}$ of $W$ ($1 \le i \le n-1$), the resulting strict transform $W^{(n)} \subset X^{(n)}$ of $W$ is regular by Proposition \ref{resbl-1}. The inverse image $Y^{(n)} \subset X^{(n)}$ of $Y$ endowed with reduced subscheme structure remains a reduced divisor with simple normal crossings by Proposition \ref{blow-upncd}. Finally applying Proposition \ref{curveanddivisor} to $(X,D,C)=(X^{(n)},Y^{(n)},W^{(n)})$, we obtain a desired proper morphism $\pi : X' \to X$. This completes the proof of Theorem \ref{res1}.
\qed

\end{document}